%% file: critii.tex
%

%
%
\magnification=\magstephalf
\input amstex
\documentstyle{amsppt}
\pagewidth{6.5truein}
\pageheight{8.9truein}
\ifx\refstyle\undefinedRLD\else\refstyle{C}\fi

\define\jalg{{{\Cal A}_j}}
\define\jalgp{{{\Cal P}_j}}
\define\elems{{\Cal E}_\lambda}
\define\Ainf{A_\infty}
\define\Binf{B_\infty}
\define\Pinf{P_\infty}
\define\freeone{{\Cal F}_{\{j\}}}
\define\crit{\operatorname{cr}}
\define\twoseq#1/#2/#3/{${#1}$~has critical sequence beginning
      ${#2}\mapsto{#3}$}
\define\threeseq#1/#2/#3/#4/{${#1}$~has critical sequence beginning
      ${#2}\mapsto{#3}\mapsto{#4}$}
\define\fourseq#1/#2/#3/#4/#5/{${#1}$~has critical sequence beginning
      ${#2}\mapsto{#3}\mapsto{#4}\mapsto{#5}$}
\define\jsub#1{j_{(#1)}}
\define\jsup#1{j^{[#1]}}
\define\restrict{\restriction}
\define\Lrestrict{\mathop{\overset*\to\cap}}

\define\ordset#1#2{S_{#1,#2}}
\define\rightblack{ \null\nobreak\hfill$\blacksquare$}
\define\rightsquare{ \null\nobreak\hfill$\square$}
\define\QED{\rightblack\enddemo}
\define\QNED{\rightblack\csname endproclaim\endcsname}
\define\QEd{\rightsquare\enddemo}
\define\procl#1{\csname proclaim\endcsname{#1}}
\define\bs{\mathbin{\bar*}}
\define\Lequiv#1/{\mathrel{\overset{#1}\to=}}
\define\bmodp{\bmod^{\!\!\prime\,\,}}
\define\subalg#1{\langle #1\rangle}
\define\kk{{\bar k}}

\define\Dehornoy{1}
\define\Dougherty{2}
\define\DougJech{3}
\define\DrapalI{4}
\define\DrapalII{5}
\define\LaverI{6}
\define\LaverII{7}
\define\SoloReinKana{8}
\define\Wehrung{9}

\topmatter
\title Critical points in an algebra of elementary embeddings, II\endtitle
\author Randall Dougherty\endauthor
\affil Ohio State University\endaffil
\date March 7, 1995 \enddate
\thanks The author was supported by NSF grant number DMS-9158092 and
by a grant from the Sloan Foundation.\endthanks
\address Department of Mathematics, Ohio State University,
Columbus, OH 43210\endaddress
\abstract
In this paper, we continue the study of a left-distributive algebra
of elementary embeddings from the collection of sets of rank less
than~$\lambda$ to itself, as well as related finite left-distributive
algebras (which can be defined without reference to large cardinals).
In particular, we look at the critical points (least ordinals moved)
of the elementary embeddings; simple statements about these ordinals
can be reformulated as purely algebraic statements concerning the left
distributive law.  Previously, lower bounds on the number of critical
points have been used to show that certain such algebraic statements,
known to follow from large cardinals, require more than Primitive
Recursive Arithmetic to prove.  Here we present the first few steps
of a program that, if it can be carried to completion, should give
exact computations of the number of critical points, thereby showing
that hypotheses only slightly beyond Primitive Recursive Arithmetic
would suffice to prove the aforementioned algebraic statements.
\endabstract
\endtopmatter

\document

\input critii1

\input critii2

\enddocument

%% file: critii1.tex
%

\head 1. Introduction \endhead
Let $V_\lambda$ be the collection of sets of rank less than $\lambda$, where
$\lambda$ is some fixed limit ordinal.  The assumption that there is a
nontrivial elementary embedding from $V_\lambda$ to $V_\lambda$ is an
extremely strong large cardinal hypothesis~\cite{\SoloReinKana}.  But
once one such embedding is known to exist, more of them can be obtained
by applying embeddings to each other.  If $a$ and $b$ are two such
embeddings, then one cannot literally apply $a$ to $b$ since $b \notin
V_\lambda$, but one can apply $a$ to initial segments of $b$, so we
define $a(b)$ to be $\bigcup_{\alpha < \lambda} a(b\restrict
V_\alpha)$.  Then $a(b)$ will also be an elementary embedding from
$V_\lambda$ to $V_\lambda$.

Let $\elems$ be the set of all nontrivial elementary embeddings from
$V_\lambda$ to $V_\lambda$, equipped with the binary operation of
application.  Elementarity implies that this operation is left
distributive over itself: $a(b(c)) = a(b)(a(c))$.  Assuming
$\elems$ is nonempty, let $j$ be a fixed nontrivial elementary embedding from
$V_\lambda$ to $V_\lambda$,  and let $\jalg$ be the subalgebra of
$\elems$ generated by $j$.

One may also consider the operation of composition of embeddings,
in which case one can list four laws satisfied by the two operations:
$$(a\circ b)\circ c = a\circ(b\circ c),\quad (a\circ b)(c) = a(b(c)),
\quad a(b\circ c) = a(b) \circ a(c), \quad a\circ b = a(b) \circ a.
\tag 1.1$$  Let $\jalgp$ be the subalgebra of $\elems$ generated
by the fixed embedding $j$ using both application and composition.

Laver~\cite{\LaverI} has shown that $\jalg$ is a free
algebra on the generator $j$ with respect to the left distributive law,
and $\jalgp$ is free with respect to \thetag{1.1}.
He also showed that $\jalg$ has an algebraic property called
irreflexivity; this gives a proof that the free left-distributive
algebra on one generator is irreflexive, assuming an extreme
large cardinal hypothesis.  The question of whether the large
cardinal is necessary for this purely algebraic result was resolved
by Dehornoy~\cite{\Dehornoy}, who proved the result without using any
large cardinals; in fact, Dehornoy's proof goes through in a very
weak system of arithmetic (Primitive Recursive Arithmetic).

Any embedding $k \in \jalg$ maps ordinals to ordinals in a strictly
increasing manner, so $k(\alpha) \ge \alpha$ for all $\alpha <
\lambda$; the least $\alpha$ such that $k(\alpha) > \alpha$ is called
the {\it critical point} of $k$, and denoted by $\crit k$.  This
$\alpha$ is inaccessible, measurable, etc., and $k(x) = x$ for $x \in
V_\alpha$~\cite{\SoloReinKana}.  Elementarity implies that $k'(\crit k)
= \crit{(k'(k))}$ and $k'(k(\beta)) = k'(k)(k'(\beta))$; using these
rules, we can obtain many ordinals below $\lambda$ as critical points
of members of $\jalg$.  To start with, let $\kappa_0 = \crit j$ and
$\kappa_{n+1} = j(\kappa_n)$ for $n \in \omega$; the ordinals
$\kappa_n$ form a strictly increasing sequence (called the {\it
critical sequence} of $j$), and all of them are critical points of
members of $\jalg$ (if $\kappa_n = \crit k$, then $\kappa_{n+1} = \crit
{j(k)}$).  But these are not the only ordinals which
occur as critical points.  In fact, if $F(n)$ is the number of critical
points of members of $\jalg$ lying below $\kappa_n$, then $F$ grows
faster than any primitive recursive function~\cite{\Dougherty}.  On
the other hand, Laver and Steel showed that $F(n)$ is finite~\cite{\LaverII}.

Laver~\cite{\LaverII} has also shown that the algebra $\jalg$ has as 
homomorphic images a sequence of {\sl finite\/} left-distributive algebras
$A_n$ which can be defined by explicit construction, without using
strong hypotheses; they can also be defined by purely algebraic
methods.  Then the results about $\jalg$ (in particular, the Laver-Steel
theorem) have consequences which are purely algebraic statements
about these finite structures; it is again natural to ask whether
large cardinals are necessary for the proof of these statements.
This question is currently open, but it is known that at least
some strength is needed: the statements cannot be proved in
Primitive Recursive Arithmetic~\cite{\DougJech}.

The statements about critical points mentioned above have natural
translations into the context of the finite algebras $A_n$; in particular,
one can define the function $F(n)$ in this context.  Then the statement
that $F(n)$, as computed from the finite algebras, is finite for all $n$
can be viewed as the algebraic version of the Laver-Steel theorem.
A more natural statement equivalent to this one is that every
term built up from one generator using the left-distributive
operation will turn out to be nonzero when evaluated in some $A_n$.
An interesting particular case: the statement that $F(4)$ is finite is
equivalent to ``There is an $n$ such that $1*16\ne 0$ in $A_n$.'' One
can easily write a computer program to search for such an $n$;
this program is known to terminate, but the only known proof
at present requires extreme large cardinals (for this particular
result, a $4$\snug-huge cardinal will suffice).  One can also give
algebraic characterizations of the algebras $A_n$ which
allow one to translate ``$F(n)$ is finite'' into purely
algebraic statements of the form ``under the left distributive law,
one type of equation does not necessarily imply another type of equation.''

The proofs of the lower bounds on $F(n)$~\cite{\Dougherty} work
by producing infinite (or at least very large) two-dimensional
arrays of embeddings, in which a certain number of embeddings in
one column generate a much larger number of embeddings in the
preceding column, and the embeddings in column~$1$ produce
a large number of critical points.  A uniform construction of this
sort produces enough critical points to show that $F$ grows faster
than any primitive recursive function.  There is no reason
to believe that this construction produces all of the critical points;
in fact, it produces $256$ critical points below $\kappa_4$, while
a much more irregular construction gives many more.  But even the
latter construction probably does not produce all critical points
in the relevant range.

The goal of this paper is to give the first steps toward an optimal
construction of this type; that is, one which will produce {\sl all}
of the critical points of members of $\jalg$.  If such a construction
can be completed explicitly, and if the critical points $\kappa_n$ appear
at known places in the array, then one will be able to give explicit
recursive formulas for computing the function~$F$, and hence
give a constructive proof of the Laver-Steel theorem.  Furthermore,
if the methods used are simple enough, then they will be translatable
into the context of the finite algebras~$A_n$, thus showing that
the algebraic statements above can be proved in a theory much weaker
than the large cardinal assertions currently needed.  In fact, one
will probably only need to assert that the recursion used to
compute~$F$ produces a total function, and this might require
no more than Primitive Recursive Arithmetic together with the
statement ``the Ackermann function is total.''  (But if $F$ turns out
to be even faster-growing than currently known, then a somewhat stronger
axiom might be needed.)

Along the way, we will produce new structural results about the
finite algebras~$A_n$, as well as additional specific results
about the lower critical points (for instance, the $256$ critical
points below~$\kappa_4$ mentioned above turn out to be the {\sl first}
$256$ critical points).  An interesting point is that some of the
main results here were discovered through examination of
bit patterns from printouts of extensive computations in the
algebras~$A_n$; this may be the first serious example of
computer-aided research in the theory of large cardinals.
Conversely, the new results led to improved algorithms for
computing in these algebras, thus allowing the development of
software for carrying out experiments at a much higher level.

The results in this paper have a large overlap with results of
A.~Dr\'apal~\cite{\DrapalI, \DrapalII}, which were obtained independently
at about the same time, but from a quite different point of view.

\head 2. Elementary embeddings and finite algebras \endhead
In this section, we review the relevant definitions and fundamental
results concerning the elementary embedding algebra $\jalg$, the finite
algebras~$A_n$ ($n \in \omega$), and the connections between them.
Most of the results here come from Laver~\cite{\LaverI,\LaverII} and
Dougherty and Jech~\cite{\DougJech}.

\subhead Elementary embeddings and critical points \endsubhead
Section~1 gave a number of the basic definitions on this topic
(the algebra~$\elems$ and its subalgebra~$\jalg$ generated by
a fixed single embedding~$j$, the critical point $\crit k$ of an
embedding $k$, etc.).  We will usually abbreviate application
by writing just~$ki$ instead of~$k(i)$.  Longer concatenations
should be associated from the left:  $k_1k_2k_3k_4$ is short
for $((k_1k_2)k_3)k_4$, or, using parentheses for application,
$k_1(k_2)(k_3)(k_4)$.  The differing uses of parentheses for
application and for grouping (not to mention for parenthetical remarks)
do not conflict with each other, but it is convenient to remove
as many such parentheses as possible.

For any embedding $k \in \jalg$, one can construct an infinite sequence
of ordinals by starting with the critical point $\crit k$ and
repeatedly applying the embedding $k$; this sequence is called the {\it
critical sequence} of~$k$.  Since $k(\crit k) > \crit k$ and $k$ is
strictly increasing on ordinals, the critical sequence is strictly
increasing; since $k(\crit i) = \crit{(ki)}$, all ordinals in the
critical sequence are critical points of members of~$\jalg$.  By a
theorem of Kunen (in the form given by Solovay, Reinhardt,
and~Kanamori~\cite{\SoloReinKana}), the critical sequence of $k$ is
cofinal in $\lambda$.  Hence, $\lambda$ is not a `large cardinal' in
the strongest sense (it is not even regular, but has cofinality~$\omega$),
but it is a cardinal, and is a limit of critical points, which are
extremely large cardinals.

The set of critical points of members of $\jalg$ is countably infinite,
since $\jalg$ is countable and any critical sequence gives infinitely
many critical points.  (One can show that $\jalg$ and $\jalgp$ give
the same set of critical points, since every element of $\jalgp$ can
be approximated by an element of $\jalg$~\cite{\LaverII}.)
Let $\gamma_0,\gamma_1,\gamma_2,\dotsc$
be a list of the first $\omega$ members of this set, in increasing
order.  Then the Laver-Steel theorem states that all critical
points of members of $\jalg$ are in this list.  This is
a consequence of two separate results.  The first, due to
Laver~\cite{\LaverII}, comes from an analysis of the particular
sequence of embeddings $\jsub n$ defined by $\jsub 1 = j$
and $\jsub{n+1} = \jsub n j$ (that is, $\jsub n = jj\dots j$ with
$n$ $j$\snug's); it follows from results given later in this review.

\proclaim{Theorem 2.1 \rm(Laver)} For any $n > 0$,
$\crit{\jsub n} = \gamma_m$, where $m$ is the largest natural number
such that $2^m$ divides $n$. \endproclaim

The second result was a by-product of inner model research of Steel.
The original proof (see Laver~\cite{\LaverII}) is in terms of
extenders; a simplified version is given below.  When applied with
$k_0 = k_1 = \dots = j$, the result states that $\{\crit \jsub n\colon
n \ge 1\}$ is cofinal in $\lambda$; this, together with Theorem~2.1,
implies that no critical point can lie above all of the $\gamma_n$\snug's,
so the $\gamma_n$\snug's are all of the critical points.

\proclaim{Theorem 2.2 \rm(Steel)} For any $k_0,k_1,\ldots\in\elems$,
$\sup_{n \in \omega} \crit(k_0k_1\dotsm k_n) = \lambda$. \endproclaim

\demo{Proof} For $i \in \elems$ and $\gamma < \lambda$, let
$\ordset i\gamma$ be the following set of ordinals:
$$\ordset i\gamma = \{i(f)(s)\colon f,s \in
V_\gamma,\ f\text{ a function whose values are ordinals}\}.$$

\procl{Claim} If $\gamma < \lambda$ is an inaccessible cardinal,
$k,i \in \elems$, and $\crit i < \gamma$, then the order type of
$\ordset{i(k)}\gamma$ is less than the order type of $\ordset i\gamma$.
\endproclaim

\demo{Proof}  Since $\crit i < \gamma$, there must be an ordinal
$\alpha < \gamma$ such that $i(\alpha) \ge \gamma$.  Let $G$ be the function
defined on pairs $\langle f,s \rangle \in V_\alpha \times V_\alpha$ with $f$ a
function whose values are ordinals and $s \in \text{domain}(k(f))$, such that
$G(\langle f,s\rangle) = k(f)(s)$; then $\ordset k\alpha$ is the range of
$G$.  We have $|\ordset k\alpha| \le |V_\alpha| < \gamma$, so there is an
order-preserving function $H: \ordset k\alpha \to \beta$ for some $\beta <
\gamma$.  Hence, $i(H)\colon \ordset {i(k)}{i(\alpha)} \to i(\beta)$ is also
order-preserving.  Note that $\ordset {i(k)}\gamma \subseteq
\ordset {i(k)}{i(\alpha)}$.  If $\delta \in \ordset {i(k)}\gamma$, then
$\delta = i(k)(f)(s)$ for some $f,s \in V_\gamma$ with $f$ a
function whose values are ordinals; then
$$\align i(H)(\delta) &= i(H)(i(k)(f)(s))\\
&= i(H)(i(G)(\langle f,s\rangle)) \\
&= i(H \circ G)(\langle f,s\rangle). \endalign$$
Since $H \circ G$ is in $V_\gamma$, $i(H)(\delta)$ is in $\ordset i\gamma$.
Therefore, $i(H)$ gives an order-preserving map from
$\ordset {i(k)}\gamma$ into $i(\beta) \cap \ordset i\gamma$, which is a
proper initial segment of $\ordset i\gamma$ [since $i(\beta) =
i(\{\langle 0,\beta \rangle\})(0) \in \ordset i\gamma$], so the order type of
$\ordset{i(k)}\gamma$ is less than that of $\ordset i\gamma$.
\QEd

Now let $i_n = k_0k_1\dots k_n$, and suppose $\sup_{n \in \omega}
\crit{i_n} < \lambda$.  Then, by Kunen's theorem,
there is an inaccessible $\gamma < \lambda$ such
that $\sup_{n \in \omega} \crit{i_n} < \gamma$.  Let $\theta_n$ be the order
type of $\ordset{i_n}\gamma$; then the claim gives $\theta_0 > \theta_1 >
\theta_2 > \ldots$, contradiction. \QED

This, along with Kunen's theorem (which can be viewed as a consequence
of the Laver-Steel theorem, but was needed for the proof of
that result), is almost the only place in the
development here where any real set theory is used; the other arguments
used could be performed in an abstract context where we have objects
called `embeddings,' other objects called `ordinals,' and some
operations and relations connecting them which satisfy a few basic
axioms~\cite{\DougJech}.

Some of the more basic facts about embeddings and ordinals are:
each embedding gives a strictly increasing mapping from ordinals
to ordinals, and $k(\alpha) \ge \alpha$ for all $\alpha$;
$\crit k$ is the least ordinal moved by $k$;
$k(i(\alpha)) = ki(k(\alpha))$ and $k(\crit i) = \crit(ki)$.
These facts can be applied to the initial critical points~$\gamma_n$
without using the Laver-Steel theorem.  This gives the following
results, which are stated in a somewhat ugly form for reasons
which will be given later.

\proclaim{Proposition 2.3} Let $k,i \in \jalg$ and $n,m \in \omega$.  Then:
\roster
\item"$\bullet$" If $k(\gamma_n) \ge \gamma_m$, then $k(\gamma_{n+1})
   \ge \gamma_{m+1}$.
\item"$\bullet$" If $\crit i \ge \gamma_n$ and $k(\gamma_n) \ge \gamma_m$,
   then $\crit{(ki)} \ge \gamma_m$.
\item"$\bullet$" If $\crit i \not\ge \gamma_{n+1}$ and $k(\gamma_n)
   \not\ge \gamma_{m+1}$, then $\crit{(ki)} \not\ge \gamma_{m+1}$.
\item"$\bullet$" $\crit k \ge \gamma_{m+1}$ if and only if $k(\gamma_m)
   \not\ge \gamma_{m+1}$. \rightblack
\endroster \endproclaim

If we have embeddings $k,i \in \elems$, and we are only interested
in determining a limited part of $ki$, then clearly we only need a limited
part of $i$: $ki \restrict V_{k(\beta)} = k(i \restrict V_\beta)$.
This only requires a limited part of $k$ as well, namely $k \restrict V_\alpha$
for any $\alpha$ greater than the rank of $i \restrict V_\beta$.  But these
limited parts are of varying sizes; in general, it is not possible
to determine $ki \restrict V_\beta$ from $k \restrict V_\beta$
and $i \restrict V_\beta$.

Laver~\cite{\LaverI} has defined a variant form of `restriction' that
gets around this problem, by suitably restricting both the domain and
the range.  Let $k \Lrestrict V_\beta = \{(x,y) \in V_\beta \times
V_\beta \colon y \in k(x)\}$.  We again have $ki \Lrestrict
V_{k(\beta)} = k(i \Lrestrict V_\beta)$, but now one can also show
that, for any limit ordinal $\beta$, $ki \Lrestrict V_\beta$
and $(k \circ i) \Lrestrict V_\beta$ depend
only on $k \Lrestrict V_\beta$ and $i \Lrestrict V_\beta$.
We can also define the corresponding version of `agreement up
to~$\beta$': say that $k \Lequiv\beta/ k'$ if $k \Lrestrict V_\beta = k'
\Lrestrict V_\beta$.  This means that $k$ and $k'$ agree about what
happens in $V_\beta$; in particular, if $k(x) \in V_\beta$ and
$k \Lequiv\beta/ k'$, then $k'(x) = k(x)$.  Also, if $\crit i \ge \beta$,
then $i$ does not move anything in $V_\beta$, so we get
$k \Lequiv\beta/ ik$ for any $k$.  Applying these
general facts to the specific case of the ordinals $\gamma_n$ gives:

\proclaim{Proposition 2.4} Let $k,k',i \in \jalg$ and $n,m \in \omega$.  Then:
\roster
\item"$\bullet$" $\Lequiv\gamma_n/$ is an equivalence relation.
\item"$\bullet$" If $k \Lequiv\gamma_n/ k'$, then $k(\gamma_m) \ge \gamma_n$
   iff $k'(\gamma_m) \ge \gamma_n$.
\item"$\bullet$" If $k \Lequiv\gamma_n/ k'$, then $ki \Lequiv\gamma_n/ k'i$,
   $k\circ i \Lequiv\gamma_n/ k'\circ i$, and
   $i\circ k \Lequiv\gamma_n/ i\circ k'$.
\item"$\bullet$" If $k \Lequiv\gamma_n/ k'$ and $i(\gamma_n) \ge
   \gamma_m$, then $ik \Lequiv\gamma_m/ ik'$. 
\item"$\bullet$" If $\crit i \ge \gamma_n$, then $k \Lequiv\gamma_n/ ik$.
   \rightblack
\endroster \endproclaim

\subhead The finite left-distributive algebras $A'_n$ and $A_n$ \endsubhead
This subsection is a summary of results from sections 3 and~4 of Dougherty
and~Jech~\cite{\DougJech}; see that paper for full proofs.

For each natural number~$n$, the
algebra $A'_n$ is defined to be the set
$\{1,2,\dots,2^n\}$ with a binary operation
$*_n$ satisfying the following recursive formulas:
$$ 2^n *_n b = b;\tag2.1a$$
if $a<2^n$, then
$$a *_n 1 = a+1;\tag2.1b$$
if $a<2^n$ and $b<2^n$, then
$$a *_n (b+1) = (a*_nb)*_n(a+1).\tag2.1c$$
This is a valid recursion, because the following inductive condition is
maintained:
$$a*_nb > a \qquad\text{if $a<2^n$.}\tag2.2\snug${}'$ $$

We will want to use reduction modulo $2^n$ to map natural numbers into
the set $\{1,2,\dots,2^n\}$, rather than the usual
$\{0,1,\dots,2^n-1\}$; let $x \bmodp 2^n$ be the unique
member of $\{1,2,\dots,2^n\}$ which is congruent to~$x$ modulo~$2^n$.
In particular, $0 \bmodp 2^n$ will be~$2^n$.

For any fixed $a$ in $A'_n$, consider the sequence $a *_n1,
a*_n2,\dots,a*_n2^n$ in~$A'_n$.  If $a=2^n$, this sequence is just
$1,2,\dots,2^n$.  If $a<2^n$, then the sequence begins with~$a+1$, and
(by~(2.1c)) each member is obtained from its predecessor by operating
on the right by~$a+1$; hence, by~\thetag{2.2\snug${}'$}, the sequence must be
strictly increasing as long as its members remain below~$2^n$.  Once
$2^n$ is reached (as must happen in at most $2^n-a$ steps), the next
member will be $a+1$ again, and the sequence repeats.  Therefore, the
sequence $a *_n1, a*_n2,\dots,a*_n2^n$ is periodic (as long as it
lasts); each period is strictly increasing from~$a+1$ to~$2^n$.  We
will refer to the number of terms in each period of this sequence as
{\it the period length of~$a$ in~$A'_n$}, or
{\it the period of~$a$ in~$A'_n$}.  (The period of~$2^n$ in~$A'_n$
is~$2^n$.)

One can prove the following results by simultaneous induction on $n$:

\proclaim{Proposition 2.5 \rm\cite{\DougJech}}
\roster
\item"{\rm(a)}" The period of any~$a$ in~$A'_n$ is a power
of~$2$; equivalently, $a*_n2^n = 2^n$ for all~$a$.

\item"{\rm(b)}" The formulas \thetag{2.1} hold modulo\/~$2^n$ in $A'_n$
even when $a$ or~$b$ is $2^n$.

\item"{\rm(c)}" Reduction modulo $2^n$ is a homomorphism from $A'_{n+1}$ to~$A'_n$:
$$(a *_{n+1} b) \bmodp 2^n = (a
\bmodp 2^n) *_n (b \bmodp 2^n)$$ for all $a,b$ in $A'_{n+1}$.

\item"{\rm(d)}" For any $a<2^n$ in $A'_n$, if $p$ is the period of $a$
in $A'_n$, then the period of~$a+2^n$ in~$A'_{n+1}$ is also $p$, and
the period of~$a$ in~$A'_{n+1}$ is either $p$ or $2p$. The period of
$2^n$ in $A'_{n+1}$ is $2^n$. \rightblack
\endroster\endproclaim

Given these properties of $A'_n$, the proof that the left distributive
law holds in $A'_n$ is a straightforward triple induction.

Let $a$ be an element of $A'_n$, with period $2^m$.  Then the left
distributive law implies that the map $b \mapsto a*_nb$ is a
homomorphism from $A'_n$ to $A'_n$, with range
$S_a = \{a*_n1,a*_n2,\dots,a*_n2^m\}$, so $S_a$ is a subalgebra of $A'_n$.
Since $A'_n$ is generated by $1$, $S_a$ is generated by $a*_n1$.
Furthermore, $a*_nb = a*_nb'$ if and only if $b \equiv b' \pmod {2^m}$,
so $S_a$ is isomorphic to the algebra $A'_n$ with elements identified
if they are congruent modulo~$2^m$; but this is just $A'_m$, by
repeated application of Proposition~2.5(c).  Since any element
$\bar a$ of $A'_n$ is of the form $a*_n1$ for some $a$
(namely $a = \bar a - 1$, or $2^n$ if $\bar a = 1$), we see
that any element of $A'_n$ generates a subalgebra which is isomorphic
to $A'_m$ for some $m$ (and the isomorphism is just the unique
order-preserving bijection between the subalgebra and $A'_m$).

One can define an algebra $A_n$ which looks just like $A'_n$, except that
$2^n$ is renamed $0$.  
This has no effect on the structure of the algebras, but it affects statements
referring to the ordering of the elements of the algebra.  In particular,
(2.2\snug${}'$) becomes:
$$\text{either}\quad a*_nb = 0 \quad \text{or} \quad a*_nb > a.\tag2.2$$
Also, the ordinary $\bmod$ operation now gives the homomorphism from
$A_{n+1}$ to $A_n$.  The element $0$ satisfies the identities
$a*_n0 = 0$ and $0*_na = a$.

The algebras $A_n$ (or $A'_n$) have a purely algebraic characterization:
$A_n$ is a free algebra on one generator $1$ subject to the left
distributive law and the equation $1 = 1*1*\dots*1$, where the product
has $2^n+1$ $1$\snug's, associated from the left.  This is proved in
Wehrung~\cite{\Wehrung}, with an error that is corrected in
Dr\'apal~\cite{\DrapalI}.

Since we have homomorphisms from $A_{n+1}$ to $A_n$ for all $n$,
we can construct an inverse limit of the algebras $A_n$; let
$\Ainf$ be the subalgebra of this inverse limit generated by the
element $1$.  Alternatively, one can define $\Ainf$ to be the
set of equivalence classes of words in the generator $1$ using the
operation $*$, where two words are considered equivalent if they
evaluate to the same result in $A_n$ for all $n$.  Then $\Ainf$ is
a countably infinite left-distributive algebra on one generator.
One can show that $\Ainf$ is free if and only if every word in
the generator $1$ evaluates to a nonzero result in $A_n$ for some $n$.

One can define a second operation $\circ_n$ on $A_n$ or $A'_n$
by the formula $$a \circ_n b = (a *_n (b+1)) - 1, \tag 2.3$$
where the addition and subtraction are performed modulo $2^n$.
Straightforward inductions show that the resulting two-operation
algebras, called $P_n$ and $P'_n$, satisfy Laver's laws \thetag{1.1}.
One also has the same homomorphisms as before from $P_{n+1}$ to
$P_n$, so one can take an inverse limit and then let $\Pinf$ be
the subalgebra generated by $1$; in contrast to the finite algebras,
the infinite algebra $\Pinf$ has elements that are not in $\Ainf$
(assuming $\Ainf$ is free).

\subhead Connections between finite algebras and embeddings \endsubhead
The review so far has been rather distorted historically.  In fact,
the finite algebras $A_n$ were not discovered until Laver proved that
certain homomorphic images of the elementary embedding algebra
$\jalg$ are finite, and determined their structure explicitly.
This result can now be stated as follows.  Recall that, for any
limit ordinal $\beta < \lambda$, $ki\Lrestrict V_\beta$ and
$(k \circ i)\Lrestrict V_\beta$ depend
only on $k\Lrestrict V_\beta$ and $i\Lrestrict V_\beta$, so the
restrictions $k\Lrestrict V_\beta$ for $k \in \jalg$ form an algebra
which is a homomorphic image of $\jalg$, and hence is
left-distributive.

\proclaim{Theorem 2.6\/ \rm(Laver~\cite{\LaverII})}
For any $n$, the algebra $\{k\Lrestrict V_{\gamma_n}\colon k \in \jalg\}$
is finite, of size $2^n$, and is isomorphic to $A'_n$,
with $j\Lrestrict V_{\gamma_n}$ corresponding to\/~$1$.
Similarly, $\{k\Lrestrict V_{\gamma_n}\colon k \in \jalgp\}$
has size $2^n$ and is isomorphic to $P'_n$,
with $j\Lrestrict V_{\gamma_n}$ corresponding to\/~$1$. \QNED

So there is a natural homomorphism from $\jalg$ to $A'_n$ which sends
$j$ to $1$.  In analogy with the natural homomorphisms from $A'_m$ ($m>n$)
to $A'_n$, we will use the notation $k \mapsto k \bmodp 2^n$ for this
new homomorphism.  When we are working with $A_n$ instead of $A'_n$,
the homomorphism is $k \mapsto k \bmod 2^n$; this again sends $j$ to
the generator $1$ (except for $n=0$, in which case it sends $j$ to~$0$,
of course).  The `identity' restricted embedding, namely $k \Lrestrict
V_{\gamma_n}$ where $\crit k \ge \gamma_n$, corresponds to $2^n$ in
$A'_n$ or $0$ in $A_n$.

The above homomorphisms give a natural homomorphism from $\jalg$ onto $\Ainf$.
Since $\jalg$ is known to be free, we have that $\Ainf$ is free
iff this homomorphism is injective.

It turns out that, even without assuming the large cardinals necessary
to define $\jalg$, one can talk about critical points of `embeddings'
purely in the context of the algebras $A_n$.  Instead of
the elementary embedding algebra $\jalg$, we can just use the free
left-distributive algebra $\freeone$ on one generator $j$.
For any $k \in \freeone$, we can express $k$ as a word in the
generator $j$; define $k \bmodp 2^n$ to be the result of evaluating
this word in $A'_n$ with $j$ set to $1$.  (Since $A'_n$ is left-distributive,
two words giving the same $k \in \freeone$ will give the same result
for $k \bmodp 2^n$.)  Define $k \bmod 2^n$ analogously.  These maps
are homomorphisms from $\freeone$ to $A'_n$ and $A_n$, and they
give a homomorphism from~$\freeone$ onto $\Ainf$, which is an
isomorphism iff $\Ainf$ is free.  Note that $\jsub m \bmod 2^n
= m \bmod 2^n$, where the right side is the ordinary `$\bmod$' operation
on numbers.

Now make the following formal definitions: $\crit k \ge \gamma_n$ is
defined to mean $k \bmod 2^n = 0$; and $k(\gamma_m) \ge \gamma_n$ is
defined to mean $(k\jsub{2^m}) \bmod 2^n = 0$ (equivalently, the period
length of $k \bmod 2^n$ in~$A_n$ is at most $2^m$).  Also, say that $k
\Lequiv\gamma_n/ k'$ if $k \bmod 2^n = k' \bmod 2^n$.

We recalled earlier that $\Ainf$ is free if and only if, for every $k \in
\freeone$, there is a largest $n$ such that $k \bmod 2^n = 0$.
In this case, we naturally say that $\crit k = \gamma_n$ for this $n$.
Similarly, we say that $k(\gamma_m) = \gamma_n$ if $n$ is largest
such that $k(\gamma_m) \ge \gamma_n$.  One can now show~\cite{\DougJech}
that all of the basic laws about the critical point and application
operations (e.g., $\crit{(ki)} = k(\crit i)$) work just as well
for these new definitions.

If we do not assume that $\Ainf$ is free, then there may be elements
$k$ of $\freeone$ for which the formal statement $\crit k \ge \gamma_n$
holds for all $n$;
we might say that ``$\crit k$ lies above all of the $\gamma_n$\snug's.''
Similarly, $k(\gamma_m)$ might lie above all of the $\gamma_n$\snug's.
However, if one is careful about stating assertions so as to allow
for this possibility, then one can continue to perform the standard
manipulations with critical points to get valid results about
the algebras $A_n$, even without assuming that $\Ainf$ is free.

One way to make sure that this works is to state all such
manipulations in terms of the primitive notions $\crit k \ge \gamma_n$
and $k(\gamma_m) \ge \gamma_n$.  This is the reason that Propositions~2.3
and~2.4 were stated in a somewhat unnatural way: only these primitive
notions are used.  It is now not hard to see that the proofs in
section~7 of Dougherty and~Jech~\cite{\DougJech} work to show that
the assertions in Propositions~2.3 and~2.4 hold in the finite
algebra context, even without the assumption that $\Ainf$ is free.

One can summarize these results as follows:

\proclaim{Proposition 2.7} The following results hold both in the context
of elementary embeddings and in the context of the finite algebras $A_n$ and
$P_n$:
\roster
\item"$\bullet$" The left distributive law and the identities (1.1).
\item"$\bullet$" Propositions 2.3 and 2.4.
\item"$\bullet$" $\crit k \ge \gamma_n$ if and only if $k \bmod 2^n
   = 0$.
\item"$\bullet$" $k \Lequiv\gamma_n/ k'$ if and only if $k \bmod 2^n
   = k' \bmod 2^n$.
\item"$\bullet$" The mapping $k \mapsto k \bmodp 2^n$ is a homomorphism
   from $\jalg$ (or $\freeone$) to $A'_n$ which sends $j$ to $1$.
\endroster
In the finite algebra context, these results are provable without any
large cardinal or freeness assumptions (in fact, in Primitive Recursive
Arithmetic). \QNED

Using these facts, one can perform an argument in terms of elementary
embeddings and then easily translate it into a result about
finite algebras.

For example, consider the following result of Laver.

\proclaim{Lemma 2.8 \rm(Laver~\cite{\LaverI, Lemma 5})} For any
$k,e_0,e_1,\dots,e_n \in \elems$, we have
$$ke_0e_1\dots e_n \Lequiv\theta_n/ k(e_0e_1\dots e_n),$$ where
$\theta_n$ is the minimum of $ke_0e_1\dots e_i(\crit k)$ for
$i < n$. \endproclaim

\demo{Proof}  Induct on $n$.  We have $e_n \Lequiv\crit k/ ke_n$,
so $ke_0e_1\dots e_n \Lequiv ke_0e_1\dots e_{n-1}(\crit k)/
ke_0e_1\dots e_{n-1}(ke_n)$.  By the induction hypothesis,
$ke_0e_1\dots e_{n-1} \Lequiv\theta_{n-1}/ k(e_0e_1\dots e_{n-1})$,
so $$ke_0e_1\dots e_{n-1}(ke_n) \Lequiv\theta_{n-1}/
k(e_0e_1\dots e_{n-1})(ke_n) = k(e_0e_1\dots e_n).$$  Putting these
together gives $ke_0e_1\dots e_n \Lequiv\min(ke_0e_1\dots e_{n-1}(\crit k),
\theta_{n-1})/ k(e_0e_1\dots e_n)$, which is the desired result. \QED

This argument goes through verbatim in the finite algebra context
if $\Ainf$ is free.  If we do not assume this, then we have to allow
for the possibility that $\crit k$ and/or $\theta_i$ is above
all of the $\gamma_n$\snug's, but this causes no difficulties in the
argument.  Alternatively, one can rephrase the result as ``if
$\crit k \ge \gamma_r$ and $ke_0e_1\dots e_i(\gamma_r) \ge \gamma_m$
for all $i < n$, then $ke_0e_1\dots e_n \Lequiv\gamma_m/ k(e_0e_1\dots e_n)$''
and modify the proof accordingly.

So, for the rest of this paper, we will phrase many results and
proofs in the language of elementary embeddings, and then draw
conclusions about the finite algebras; these conclusions will be
valid without large cardinal or freeness assumptions as long
as we use only the basic methods that are justified by Proposition~2.7.
One can go in the reverse direction as well; for instance,
Proposition~2.5(d) implies that, if $k \in \jalg$ and
$k \bmod 2^{n+1} \ge 2^n$, then $k(\gamma_m) \ne \gamma_n$ for all $m$.

Since we can talk about critical points in the finite algebra
context, we can define the function~$F$ there.  In fact
$\kappa_n$ is the critical point of the embedding
$\jsup n = j(j(\dots(j(j))\dots))$, where there are $n+1$ $j$\snug's.
So $F(n)$ should be the $m$ such that $\crit{\jsup n} = \gamma_m$;
that is, the largest $m$ such that $\jsup n \bmod 2^m
= 0$.  The Laver-Steel theorem in this context is the statement
that $F(n)$ exists (as a finite number) for all $n$; that is,
for every $n$ there is an $m'$ such that $\jsup n \bmod 2^{m'} \ne 0$.
This statement is also equivalent to the assertion that $\Ainf$ is
free~\cite{\DougJech}.

Using Wehrung's characterization of $A_n$, we can get purely algebraic
versions of these statements.  For instance, ``$F(n)$ is finite''
can be restated as: There exists an $m>0$ such that the left distributive
law and the equation $\jsub {m+1} = j$ do {\sl not} imply
the equation $\jsup n j = j$.

Consider a particular $n$, say $n=4$.  If $F(4)$ is finite, then its
value can be computed; one simply has to write a program which
evaluates $j(j(j(jj)))$ for $j=1$ in the finite algebras
$A_1,A_2,\dotsc$ successively, and stops when it gets a nonzero result.
This gives a very straightforward computer program which is
known to terminate, but the only known proof of termination uses
an extreme large cardinal hypothesis.  (Since extremely large
lower bounds on $F(4)$ are known~\cite{\Dougherty}, it is
hopeless to actually
run such a program to determine $F(4)$.)

\subhead The finite algebras in backward form \endsubhead
Let $B_n$ be the algebra on
the set $\{0,1,\dots,2^n-1\}$ with one binary operation $\bs_n$
defined by the following recursive equations: $0 \bs_n x = x$,
$x \bs_n (2^n-1) = x-1$ for $x>0$, and
$x \bs_n (y-1) = (x \bs_n y) \bs_n (x-1)$ for $x,y > 0$.  Comparison
with the defining equations for $A'_n$ shows that the map
$x \mapsto 2^n-x$ is an isomorphism between $A'_n$ and $B_n$, so
results about~$A'_n$ can be translated into results about $B_n$
and vice versa.  This reversal obviously has no effect on the
structure of the algebras, but it enables one to see certain numerical
patterns that are more obscure in the original algebras.

Let us start with some direct translations of results from $A'_n$.  The
algebra $B_n$ is left-distributive and generated by the element
$2^n-1$.  If $x,y \in B_n$ and $x \ne 0$, then $x \bs_n y < x$.
Reduction modulo $2^n$ is a homomorphism from $B_N$ ($N \ge n$) to
$B_n$.  For any nonzero $x$ in $B_n$, the sequence $x \bs_n 0, x\bs_n
1, \dotsc$ is periodic with period length $2^m$ for some $m < n$; each
period of the sequence starts with $0$ and strictly increases to
$x-1$.  The set $\{x \bs_n i\colon 0 \le i < 2^m\}$ is equal to the
subalgebra generated by $x-1$, and this subalgebra is isomorphic to
$B_m$.  Clearly $x$ is the largest element of the subalgebra generated
by $x$.  If $y$ is in the subalgebra generated by $x$, then either
$y=x$ or the subalgebra generated by $y$ is strictly smaller than that
generated by $x$.

In addition to the homomorphisms above, $B_n$ is literally a subalgebra
of $B_N$ for $N \ge n$.  (To see this quickly, note that the inequality
$x \bs_N y < x$ implies that $\{0,1,\dots,2^n-1\}$ is closed under
the operation $\bs_N$; reduction modulo $2^n$ is the identity on this
set and must send $\bs_N$ to $\bs_n$, so $\bs_N$ matches $\bs_n$
on this set.)  Therefore, we can drop the subscript and just write
$\bs$ for the operation on $B_n$.

One can take the inverse limit of the algebras $B_n$, getting the same
result as for $A'_n$ in reversed form.  One can also take the union of
the algebras $B_n$ to get a countably infinite algebra $\Binf$ with
operation $\bs$ which is still left-distributive.  This algebra is
quite different from the algebra $\Ainf$; in particular, it is not
free, since the subalgebra generated by any element is finite.  Since
the period length of $x$ in $B_n$ does not depend on $n$ (it is the
least positive $l$ such that $x \bs l = 0$), the sequence $x \bs 0, x
\bs 1, \dotsc$ is still periodic in $\Binf$ with the same period.

Let $\subalg x$ denote the subalgebra of $\Binf$ (or of $B_n$ for any
$n$ with $2^n>x$) generated by $x$.  Reduction modulo $2^n$ maps
$\subalg x$ onto $\subalg{x \bmod 2^n}$.  When one moves from $B_n$ to
$B_{n+1}$, if $x \bmod 2^{n+1} = x \bmod 2^n$, then $\subalg{x \bmod
2^{n+1}} = \subalg{x \bmod 2^n}$; if $x \bmod 2^{n+1} = (x \bmod 2^n) +
2^n$, then, since reduction modulo $2^n$ maps $\subalg{x \bmod
2^{n+1}}$ (whose size is a power of $2$) onto $\subalg{x \bmod 2^n}$,
$\subalg{x \bmod 2^{n+1}}$ consists of either all or half of the
numbers $y$ and $y + 2^n$, where $y$ is in $\subalg{x \bmod 2^n}$.  In
particular, one can now prove by induction on $n$ that, if $y$ is in
$\subalg{x \bmod 2^n}$, then for any position at which the binary
representation of $y$ has a $1$, the binary representation of $x$ also
has a $1$.  In particular, this holds when $y \in \subalg x$.  It
follows that the size of $\subalg x$ (i.e., the period of $x+1$) is at
most $2^r$, where $r$ is the number of $1$ bits in the binary
representation of $x$.  Numbers $x$ for which this bound is attained
will be very useful later (for instance, because the isomorphism
from $B_r$ to $\subalg x$ is clear: one uses the $1$ bits
of $x$ as a template, and transfers the bits of a given $r$-bit
number to the specified locations, filling in $0$\snug's elsewhere).

One can translate the preceding results back into the language of
$A'_n$ or $A_n$, obtaining the following.  The map $x \mapsto x+2^n$ is
a homomorphism from $A'_n$ to $A'_{n+1}$.  For any $x$ in $A_n$, if we
write the binary expansion of $x$ with enough leading $0$\snug's to
give exactly $n$ bits in all, and $r$ of these bits are $0$, then the
period length of $x$ is at most $2^r$; in fact, for any $y$, and for
any position at which the binary representation of $x$ has a $1$, the
binary representation of $(x*y-1) \bmod 2^n$ also has a~$1$.  These
results could have been proved directly, using results such as
Proposition~2.5, instead of via the reverse algebras.

Some of these results can be translated into statements about
embeddings and critical points.  One of these will be useful
later:

\proclaim{Lemma 2.9} Let $k$ be in $\jalg$.  If $k(\gamma_m) =
\gamma_n$, then $k \bmod 2^N$ (for $N > n$) must have a\/~$0$ in
position number $n$ (numbered from the right starting with\/~$0$)
of its binary expansion. \QNED

It is easy to see that $k(\gamma_0) = \gamma_n$
where $n$ is the position of the {\sl rightmost} $0$ in the
binary expansion of $k \bmod 2^N$, if $N > n$.  We will see later
(Proposition~6.6) that
the position of the next $0$ bit of $k \bmod N$ always gives
$k(\gamma_1)$, but later values $k(\gamma_m)$ may skip over $0$ bits.

\head 3. Column 1 \endhead
Laver showed that the critical points $\gamma_0,\gamma_1,\dotsc$, which
were defined to be the first critical points of any embeddings in
$\jalg$, are actually the critical points of a specific list of
embeddings, namely the embeddings $\jsub{2^n}$, $n=0,1,\dotsc$.
He also showed~\cite{\LaverII} that $\jsub{2^n}(\gamma_n) = \gamma_{n+1}$.
We can put these two statements together by saying
``\twoseq\jsub{2^n}/\gamma_n/\gamma_{n+1}/.''  Similar notation
will be used for longer sequences: ``\threeseq k/\alpha/\alpha'/\alpha''/''
means that $\crit k = \alpha$, $k(\alpha) = \alpha'$, and
$k(\alpha') = \alpha''$.

This sort of `linked list' of embeddings and critical sequences is
a key tool in Dougherty~\cite{\Dougherty}.  The main construction in that
paper is an infinite array of embeddings, in which each column is
a list of embeddings whose critical sequences are linked in a certain
way (which varies from column to column); each column is used iteratively
to produce the preceding column.
Actually, the paper gives two constructions of this sort, one in
section~2 and the other in section~4; these will be
referred to here as `the regular construction' and
`the irregular construction.'  For instance, column 1
(in both variations of the construction) is a list of embeddings
$a_0,a_1,\dotsc$ such that \threeseq a_n/\zeta/\alpha_n/\alpha_{n+1}/.
Given such a list, one can produce $2^n$
critical points below~$\alpha_n$ \cite{\Dougherty,~Lemma~2}.

As stated in the introduction, the goal here is to produce an `optimal'
construction of this type, which yields {\sl all} of the critical
points of members of $\jalg$.  If this optimal construction is to
have a column as described above, then we must have $\zeta = \gamma_0$
and $\alpha_n = \gamma_{2^n}$ for all $n$.  The following result shows
that one can indeed produce a sequence of embeddings having this form.

\proclaim{Theorem 3.1} For any natural number $n$, the following are true:
\roster
\item"\rm(a)" For any $m < n$, there exists $\kk$ in $\jalg$
such that \threeseq\kk/\gamma_{2^n-2^{m+1}}/\gamma_{2^n-2^m}/\gamma_{2^n}/.
\item"\rm(b)" If $k \bmod 2^{2^n} \ne 0$, then $k(\gamma_{2^n}) \ge
\gamma_{2^{n+1}}$.  Furthermore, if also $k \bmod 2^{2^{n+1}} > 2^{2^n}$,
then $k(\gamma_{2^n}) > \gamma_{2^{n+1}}$.
\item"\rm(c)"
\threeseq\jsub{2^{2^n}-1}/\gamma_0/\gamma_{2^n}/\gamma_{2^{n+1}}/.
\endroster \endproclaim

\demo{Proof} By simultaneous induction on $n$.  The case $n=0$ of
(a) is vacuously true.

First, suppose (a) is true for $n$.  We prove (b) for $n$ by downward
induction on $k \bmod 2^{2^n}$.

If $k \bmod 2^{2^n} = 2^{2^n}-1$, then the rightmost $2^n$ bits of the
binary expansion of $k \bmod 2^{2^{n+1}}$ are all $1$, so at most
$2^n$ of the rightmost $2^{n+1}$ bits are $0$.  Hence, by Lemma~2.9,
the $2^n+1$ critical points $k(\gamma_0),k(\gamma_1),\dots,k(\gamma_{2^n})$
cannot all be below $\gamma_{2^{n+1}}$, so
$k(\gamma_{2^n}) \ge \gamma_{2^{n+1}}$.  If we also have
$k \bmod 2^{2^{n+1}} > 2^{2^n}$, then at most $2^n-1$ of
the rightmost $2^{n+1}$ bits of $k \bmod 2^{2^{n+1}}$ are $0$,
so we get $k(\gamma_{2^n-1}) \ge \gamma_{2^{n+1}}$ and hence
$k(\gamma_{2^n-1}) > \gamma_{2^{n+1}}$.

Now suppose $k \bmod 2^{2^n} < 2^{2^n}-1$.  Then $kj \bmod 2^{2^n} = (k
\bmod 2^{2^n})+1 > 0$, so $k(\gamma_0) = \crit{(kj)} < \gamma_{2^n}$.
On the other hand, since $k \bmod 2^{2^n} \ne 0$, we have $\crit k <
\gamma_{2^n}$, so $k(\gamma_{2^n-1}) \ge \gamma_{2^n}$.  Therefore,
there is a largest $m$ such that $k(\gamma_{2^n-2^m}) \ge
\gamma_{2^n}$, and this $m$ is less than $n$.  Apply (a) to get $\kk
\in \jalg$ such that
\threeseq\kk/\gamma_{2^n-2^{m+1}}/\gamma_{2^n-2^m}/\gamma_{2^n}/.  Then
$\crit{(k\kk)} = k(\crit \kk) = k(\gamma_{2^n-2^{m+1}})$, which is less
than $\gamma_{2^n}$ by the minimality of $m$, so $k\kk \bmod 2^{2^n}
\ne 0$.  Hence, by \thetag{2.2}, $k\kk \bmod 2^{2^n} > k \bmod 2^{2^n}$.
Therefore, we can apply the inner inductive hypothesis to get
$k\kk(\gamma_{2^n}) \ge \gamma_{2^{n+1}}$.  Since $$k(\gamma_{2^n}) =
k(\kk(\gamma_{2^n-2^m})) = k\kk(k(\gamma_{2^n-2^m})) \ge
k\kk(\gamma_{2^n}),$$ we have $k(\gamma_{2^n}) \ge \gamma_{2^{n+1}}$.
If we also have $k \bmod 2^{2^{n+1}} > 2^{2^n}$, then \thetag{2.2} gives $k\kk
\bmod 2^{2^{n+1}} > k \bmod 2^{2^{n+1}} > 2^{2^n}$, so the inner
inductive hypothesis gives $k\kk(\gamma_{2^n}) > \gamma_{2^{n+1}}$ and
hence $k(\gamma_{2^n}) > \gamma_{2^{n+1}}$.  This completes the proof
of (b) for $n$.

Next, we prove (c) for $n$ from (b) for $n$.  For any $r$ such that
$1 \le r < 2^{2^n}$, we can apply Lemma~2.8 with $k = \jsub{2^{2^n}r}$
and $e_0=e_1=\dots=j$ to get $\jsub{2^{2^n}(r+1)} \Lequiv\theta/
\jsub{2^{2^n}r}(\jsub{2^{2^n}})$, where
$\theta$ is the minimum of $\jsub{2^{2^n}r+s}(\crit k)$ for
$1 \le s < 2^{2^n}$.  Since $k \bmod 2^{2^n} = 0$, we have
$\crit k \ge \gamma_{2^n}$.  We can apply (b) to get
$\jsub{2^{2^n}r+s}(\gamma_{2^n}) > \gamma_{2^{n+1}}$ for all
of the above values of $s$, so we must have
$\theta \ge \gamma_{2^{n+1}+1}$.  Therefore,
$\jsub{2^{2^n}(r+1)} \Lequiv\gamma_{2^{n+1}+1}/
\jsub{2^{2^n}r}(\jsub{2^{2^n}})$.  Applying this for
all $r$ from $2^{2^n}-1$ down to $1$ gives
$$\align \jsub{2^{2^{n+1}}}
&\Lequiv\gamma_{2^{n+1}+1}/ \jsub{(2^{2^n}-1)2^{2^n}}\jsub{2^{2^n}} \\
&\Lequiv\gamma_{2^{n+1}+1}/
   \jsub{(2^{2^n}-2)2^{2^n}}\jsub{2^{2^n}}\jsub{2^{2^n}} \\
&\Lequiv\gamma_{2^{n+1}+1}/ \dotso \\
&\Lequiv\gamma_{2^{n+1}+1}/ \jsub{2^{2^n}}\jsub{2^{2^n}}\dots\jsub{2^{2^n}}
   \qquad\text{($2^{2^n}$ times)} \\
&\Lequiv\hphantom{\gamma_{2^{n+1}+1}}/ \jsub{2^{2^n}-1}\jsub{2^{2^n}}.
\endalign$$
So $\crit{\jsub{2^{2^{n+1}}}}$ and
$\crit{(\jsub{2^{2^n}-1}\jsub{2^{2^n}})} =
\jsub{2^{2^n}-1}(\gamma_{2^n})$ are the same if either is below
$\gamma_{2^{n+1}+1}$.  But Theorem~2.1 gives $\crit{\jsub{2^{2^n}}} =
\gamma_{2^n}$ and $\crit{\jsub{2^{2^{n+1}}}} = \gamma_{2^{n+1}}$, so we
must have $\jsub{2^{2^n}-1}(\gamma_{2^n}) = \gamma_{2^{n+1}}$.  Since
$\jsub{2^{2^n}-1}(\gamma_0) = \crit{\jsub{2^{2^n}}} = \gamma_{2^n}$,
\threeseq\jsub{2^{2^n}-1}/\gamma_0/\gamma_{2^n}/\gamma_{2^{n+1}}/, as
desired.

Finally, assuming (a) and (c) hold for $n$, we get (a) for $n+1$ as
follows.  Note that, since there are only $2^n-1$ critical points
between $\gamma_{2^n}$ and $\gamma_{2^{n+1}}$, which must include the
ordinals $\jsub{2^{2^n}-1}(\gamma_r)$ for $0 < r < 2^n$, we must have
$\jsub{2^{2^n}-1}(\gamma_r) = \gamma_{2^n+r}$ for $0 < r < 2^n$.  Now,
to get $a$ for $n+1$ and~$m$ where $m = n$, use the embedding
$\jsub{2^{2^n}-1}$.  For $m < n$, apply (a) for $n$ and $m$ to get
$\kk$; then $\jsub{2^{2^n}-1}(\kk)$ works for $n+1$ and $m$. \QED

This result makes it reasonable to believe that, in the ultimate
two-dimensional construction (if it exists), column 1 is a sequence
of embeddings $a_0,a_1,\dotsc$ such that \threeseq a_n/\gamma_0/
\gamma_{2^n}/\gamma_{2^{n+1}}/.  However, these embeddings will
probably not be precisely the embeddings $\jsub{2^{2^n}-1}$, but instead
iterative combinations of embeddings from the next column.

\head 4. A theorem about the finite algebras $A_n$ or $B_n$\endhead
This section gives a result about the left-distributive algebras $A_n$
which will be needed in order to extend the results of section~3.
It is actually stated and proved in terms of the backward algebra~$\Binf$
(the union of the finite algebras $B_n$), because the
patterns that it describes are clearer in this form.

The history of this result is interesting.  Weaker versions of it,
formulated in terms of critical points and embeddings in $\jalg$, were
conjectured as early as 1986, before the finite algebras~$A_n$ had been
discovered; but attempts to prove these results by induction were not
successful.  Much later, computer experiments displaying bit patterns
occurring in period sequences in $\Binf$ related to these conjectures
led to the statement of the stronger theorem below, which turned out to
be strong enough that an inductive proof could be pushed through.

All variables below (except $B$ and $\phi$) range over nonnegative integers.
Recall some number-theoretic notations: $2^k \mid y$ means that $y$ is a
multiple of $2^k$, and $2^k \parallel y$ means $2^k \mid y$ but
$2^{k+1} \nmid y$.

\proclaim{Theorem 4.1} Suppose $a = x + 2^{mn}y$, where $n$ is a power
of\/ $2$, $y < 2^n$, and $2^{(m+1)n} \mid x$.  Let\/ $2^l$ be the
period length of $x+1$ in $\Binf$ (i.e., the size of\/ $\subalg x$), and
suppose $l \le n$.  Then\/ $\subalg a$ is isomorphic to\/ $\subalg{
2^{l+n}-2^n+y}$, and the isomorphism can be described as follows:  for
any $i$, if\/ $(2^{l+n}-2^n + y + 1) \bs i = 2^n w + y'$ where $y' <
2^n$, then\/ $(a+1) \bs i = (x+1) \bs w + 2^{mn}y'$.  \endproclaim

The main case of interest in Theorem 4.1 (from which the full theorem
can be deduced) is the case $y=2^n-1$, which
is worth stating separately:

\proclaim{Corollary 4.2} Suppose $a = x + 2^{mn}(2^n-1)$, where $n$ is
a power of\/ $2$ and $2^{(m+1)n} \mid x$.  Let $2^l$ be the period
length of $x+1$ in $\Binf$, and suppose $l \le n$.  Then $a+1$ has
period length $2^{l+n}$, and $(a+1)\bs(2^ni+i') = (x+1)\bs i +
2^{mn}i'$ for all $i$ and all $i' < 2^n$. \endproclaim

The proof of Theorem 4.1 will be by induction on $a$.  We start by establishing
several consequences of the induction hypothesis.

\proclaim{Lemma 4.3} Assume that Theorem 4.1 (and hence Corollary 4.2) holds
for the case $a < L$, $y \ne 0$.
If $s \le s'$, $2^{s'} - 2^s < L$, and
$s' - s$ is less than or equal to some power of\/ $2$ which divides\/ $2s$,
then\/ $2^{s'} - 2^s + 1$ has period length\/ $2^{s' - s}$ (the largest
possible). \endproclaim

\demo{Proof} Induct on $s' - s$; the case $s' - s \le 1$ is trivial.
(The period of any $x > 1$ includes both $0$ and $x-1$, so its length
is at least $2$; if $s'-s = 1$, then $2^{s'-s}$ has only one $1$ bit,
so the period length of $2^{s'}-2^s+1$ is at most $2$.)
If $s' - s > 1$, let $n$ be the largest power of $2$ less than $s'-s$.
Note that $n \mid s$, so $n \mid (s+n)$; also, $s' - (s+n) \le n$.
Apply the induction hypothesis to see that $2^{s'} - 2^{s+n} + 1$ has
period length $2^{s' - (s+n)}$.  Now apply Corollary 2 with
$x = 2^{s'} - 2^{s+n}$ and $m = s/n$
to get the desired result. \QED

\proclaim{Lemma 4.4} Assume Theorem 4.1 holds for $a < L$.  If 
$x + 2^s - 1 < L$, $2^s \mid x$,
and $x+1$ has period length\/ $2^l$, where $l$ is less than or equal to
some power of\/ $2$ dividing $s$, then
$(x + 2^s)\bs(2^si + i') = (x+1)\bs i + i'$ for all $i$ and all $i'<2^s$
(so the period length of $x+2^s$ is\/ $2^{l+s}$).
\endproclaim

\demo{Proof} Induct on $s$; the case $s = 0$ is trivial.  If $s > 0$,
let $n$ be the largest power of $2$ dividing $s$; then $2n$ divides
$s -n$.  Apply Corollary 2 with $m = s/n - 1$
to see that $x + 2^s - 2^{s-n} + 1$ has period length $2^{l+n}$ and
$(x+2^s-2^{s-n}+1)\bs(2^n i + j) = (x+1)\bs i + 2^{s-n}j$
for all $i$ and all $j < 2^n$.  Now, since $l + n \le 2n$, we can apply
the induction hypothesis to get
$$\align
(x + 2^s)\bs(2^s i + 2^{s-n} j + k) &= (x + 2^s - 2^{s-n}+1)\bs(2^n i + j)+k\\
&= (x+1)\bs i + 2^{s-n}j + k
\endalign$$
for all $i$, all $j<2^n$, and all $k < 2^{s-n}$,
which gives the desired result if, given $i'$, we choose $j$ and $k$
so that $i' = 2^{s-n}j + k$. \QED

\proclaim{Lemma 4.5}
If $a = 2^nx+y$ where
$n$ is a power of\/ $2$, $x < 2^n$, and $y < 2^n-1$, then $(a+1)\bs 2^n
= 0$. \endproclaim

\demo{Proof}
This is the finite algebra form of the first sentence of Theorem~3.1(b).
%
%
\QED

\proclaim{Lemma 4.6} Assume Theorem 4.1 holds for $a < L$.  If
$a$, $x$, $y$, $n$, $m$, and $l$ are as in Theorem~4.1, $a < L$,
and $y < 2^n-1$, then $(a+1)\bs 2^n = 0$.\endproclaim

\demo{Proof} By Theorem 4.1, the period length of $a+1$ is the same as the
period length of $2^{l+n} - 2^n + y + 1$, which is at most $2^n$ by
Lemma 4.5. \QED

The preceding lemmas were consequences of Theorem 4.1 that will be useful
as applications after the theorem is proved (so that the initial
assumption becomes redundant).  The next two, however, are just parts
of the theorem itself to be used in the induction.

\proclaim{Lemma 4.7} Assume Theorem 4.1 holds for $a < L$.  If $n$ is a
power of\/ $2$, $0 < y' < 2^n$, $2^{(m+1)n} \mid x'$, $2^{(m+1)n} \mid x$,
$x'+1$ has period length at most\/ $2^n$,
and $x' + 2^{mn}y' \le L$, then $(x'+2^{mn}y')\bs(x+2^{mn}(2^n-1)) =
x'+2^{mn}(y'-1)$. \endproclaim

\demo{Proof} Apply Lemma 4.6 to $x' + 2^{mn}(y' - 1)$ to show that
$x' + 2^{mn}(y' - 1) + 1$ has period length at most $2^n$.  Then
apply Lemma 4.4 with $s = mn$ to show that
$$(x'+2^{mn}y')\bs(x+2^{mn}(2^n-1)) = (x'+2^{mn}(y'-1)+1)\bs(x/2^{mn}
+ 2^n-1).$$ But $x/2^{mn}+2^n-1 \equiv 2^n-1 \pmod{2^n}$, so
$(x'+2^{mn}(y'-1)+1)\bs(x/2^{mn}+2^n-1) = x'+2^{mn}(y'-1)$.
\QED

Define $x^{\bs(n)}$ to mean
$\overbrace{x\bs x\bs\dots\bs x}^{\text{$n$ times}}$
(associated to the left); this is
equal to $(x+1)\bs(K-n)$ if $K$ is a multiple of the period length
of $x+1$.

\proclaim{Lemma 4.8} Assume Theorem 1 holds for $a < L$.  If $n$ is a
power of\/ $2$, $2^{(m+1)n} \mid x$, $x'$ is in\/~$\subalg x$, $x' < L$,
$x \le L$, and $x+1$ has period length at most\/ $2^n$, then $$x' \bs (x
+ 2^{mn}(2^n-1)) = (x'\bs x) + 2^{mn}(2^n-1).$$ \endproclaim

\demo{Proof} Induct on $x$.
The case $x'=0$ is trivial, so assume $x' \ne 0$ and hence
$x \ne 0$.  Since $x'$ is in $\subalg x$, any $1$ bits
in the binary representation of $x'$ must correspond to $1$ bits in
the binary representation of $x$.  In particular,
this gives $x' \equiv 0 \pmod{2^{(m+1)n}}$, so
$x' \bs (x + 2^{mn}(2^n-1)) = \hat x + 2^{mn}(2^n-1)$ for
some $\hat x$ such that $\hat x < x'$ and $2^{(m+1)n} \mid \hat x$.
We need to show that $\hat x = x' \bs x$.

Since $x \ne 0$, there is a largest number $M$ such that $2^M \le x$;
clearly $x < 2^{M+1}$ and $M \ge (m+1)n$.  Now $(x \bmod 2^M)$ is a multiple of
$2^{(m+1)n}$ less than $2^M$, so $(x \bmod 2^M) + 2^{mn}(2^n-1) < 2^M
\le L$.  Since $x'$ is in $\subalg x$,
$(x' \bmod 2^M)$ is in $\subalg{x \bmod 2^M}$.
Therefore, we can apply the induction hypothesis to see that
$$(x' \bmod 2^M) \bs ((x \bmod 2^M) + 2^{mn}(2^n-1)) =
((x'\bs x) \bmod 2^M) + 2^{mn}(2^n-1).$$
Therefore, $\hat x \equiv x'\bs x \pmod{2^M}$.  Since $x'\bs x$ is a member
of $\subalg x$ and is less than $x$ (in fact,
less than $x'$), the period
length of $(x'\bs x)+1$ is less than that of $x+1$ and hence less than
$2^n$.  Since the period length of $\hat x + 1$ is at most twice
that of $(\hat x \bmod 2^M) + 1 = ((x'\bs x) \bmod 2^M) + 1$, $\hat x + 1$
has period length at most $2^n$.  Therefore, we can apply Corollary 4.2
to get $(\hat x + 2^{mn}(2^n-1))^{\bs(2^n)} = \hat x$.

Since $x' < 2^{M+1}$, $x' - 1$ has at most $M$ $1$ bits in its binary
expansion, so $x'$ has period length at most $2^M$.  We can apply Corollary 4.2
to get $((x \bmod 2^M) + 2^{mn}(2^n-1))^{\bs(2^n)} = (x \bmod 2^M)$,
so $(x + 2^{mn}(2^n-1))^{\bs(2^n)} \equiv x \pmod{2^M}$.
Finally, by the distributive law,
$$\align
\hat x &= (\hat x + 2^{mn}(2^n-1))^{\bs(2^n)} \\
&= (x' \bs (x + 2^{mn}(2^n-1)))^{\bs(2^n)} \\
&= x' \bs (x + 2^{mn}(2^n-1))^{\bs(2^n)} \\
&= x' \bs x,
\endalign$$
as desired.
\QED

\demo{Proof of Theorem 4.1} We will show that, if Theorem 4.1 holds for $a < L$,
then it holds for $a < L+1$.  Suppose $a$, $x$, $y$, $m$, $n$, and $l$ are
as hypothesized, and $a \le L$.  First, consider the case $x = 0$.
We show by induction on $z$ that $2^{mn}z \bs 2^{mn}z' =
2^{mn}(z \bs z')$ for $z,z' \le y$.  This is trivial for $z=0$.  For $z > 0$,
we have $2^{mn}z \bs 2^{mn}z' = (2^{mn}(z-1)+1)\bs z'$ by Lemma 4.4,
and $(2^{mn}(z-1)+1)\bs z' = (2^{mn}(z-1))^{\bs(2^n-z')} =
2^{mn}(z-1)^{\bs(2^n-z')} = 2^{mn}(z \bs z')$ by the induction hypothesis
on $z$.  This completes the induction on $z$.
So, in particular, $\subalg y$ is isomorphic to
$\subalg {2^{mn}y}$ via the map $z \mapsto 2^{mn}z$,
which is the desired result.

Second, we handle the case where both $x$ and $y$ are nonzero.
This implies $x < L$ and $2^{l+n} - 2^n < L$, so we can apply Lemma 4.3
to see that $2^{l+n}-2^n+1$ has maximal period length $2^l$; that is,
$2^n(i+1) \bs (2^{l+n}-2^n) = 2^ni$ for all $i < 2^l-1$.
Hence, $\subalg x$ and $\subalg{2^{l+n}-2^n}$ are
isomorphic, with the isomorphism $\phi$ mapping $(x+1)\bs w$ 
to $2^nw$ for $w < 2^l$.
This implies that $(x+1)\bs w + 1$ and $2^nw + 1$ have the same
period length for each such $w$.

The conclusion of Theorem 4.1 certainly holds for infinitely many $i$, namely
for all common multiples of the period lengths of $a+1$ and
$2^{l+n}-2^l+y+1$ (since this gives $(a+1)\bs i = 0 = (2^{l+n}-2^l+y+1)\bs i$).
Therefore, it will suffice to show that, if the conclusion holds for
$i+1$, then it holds for $i$.

Suppose the conclusion holds for $i+1$.  If $(2^{l+n}-2^n+y+1)\bs(i+1) = 0$,
then $(a+1)\bs(i+1)=0$, and hence $(2^{l+n}-2^n+y+1)\bs i = 2^{l+n}-2^n+y$
and $(a+1)\bs i = a$, so the conclusion holds (with $w = 2^l-1$ and $y'=y$).
So suppose $(2^{l+n}-2^n+y+1)\bs(i+1) = 2^nw+y' \ne 0$; then
$(a+1)\bs(i+1) = x' + 2^{mn}y'$, where $x' = (x+1)\bs w$.

Write $2^nw+y'-1 = (2^nw+y')\bs(2^{l+n}-1)$ as $2^n\hat w+\hat y$ with
$\hat y < 2^n$.  If $y' > 0$, then $\hat w = w$ and $\hat y = y'-1$;
if $y'=0$, then $\hat w = w-1$ and $\hat y = 2^n-1$.  
Now Lemma 4.7 (if $y' > 0$) or Lemma 4.8 (if $y'=0$) gives
$(x' + 2^{mn}y')\bs(x + 2^{mn}(2^n-1)) = \hat x + 2^{mn}\hat y$,
where $\hat x = (x+1)\bs\hat w$.

As noted before, the isomorphism $\phi$ ensures that
$(x+1)\bs\hat w+1$ and $2^n\hat w+1$ have the same period
length; call this length $2^{\hat l}$.  Clearly $\hat l \le l$.
We can now apply the main induction hypothesis to both
$\hat x + 2^{mn}\hat y$ and $2^n\hat w + \hat y$ to see that they both
generate subalgebras isomorphic to $\subalg{2^{\hat l + n}
- 2^n + \hat y}$ and hence
to each other, with the isomorphism $\hat \phi$
mapping
$(\hat x + 1)\bs v + 2^{mn}z$ to $((2^n\hat w + 1) \bs v) + z$ for all
valid values of $v$ and $z < 2^n$ (i.e., for all $v$ and $z$ such that
$2^nv+z$ is in $\subalg{2^{\hat l+n}-2^n+\hat y}$).

We have $2^{l+n}-2^n+y = (2^{l+n}-1)^{\bs(2^n-y)}$.  We next show that
$x+2^{mn}y$ is congruent to $(x+2^{mn}(2^n-1))^{\bs(2^n-y)}$ modulo
the period length of $x'+2^{mn}y'$.  Since $x \ne 0$, let $2^M$ be
the largest power of $2$ less than or equal to $x$.  Then we can
apply Corollary 4.2 to $(x \bmod 2^M) + 2^{mn}(2^n-1)$ to get
$$((x \bmod 2^M)+2^{mn}(2^n-1))^{\bs(2^n-y)} =
(x \bmod 2^M)+2^{mn}y.$$  Since $x'+2^{mn}y'-1$ is less than $2^{M+1}-1$,
it has at most $M$ $1$ bits in its binary expansion, so the period
length of $x'+2^{mn}y'$ is at most $2^M$, so we have the desired congruence.

Now we can write
$$ \align
(2^{l+n}-2^n+y+1)\bs i &= (2^nw'+y')\bs(2^{l+n}-2^l+y) \\
&= (2^nw'+y')\bs(2^{l+n}-1)^{\bs(2^n-y)} \\
&= ((2^nw'+y')\bs(2^{l+n}-1))^{\bs(2^n-y)} \\
&= (2^n\hat w+\hat y)^{\bs(2^n-y)},
\endalign$$
$$\align
(a+1)\bs i &= (x' + 2^{mn}y') \bs (x + 2^{mn}y) \\
&= (x' + 2^{mn}y') \bs (x+2^{mn}(2^n-1))^{\bs(2^n-y)} \\
&= ((x' + 2^{mn}y') \bs (x+2^{mn}(2^n-1)))^{\bs(2^n-y)} \\
&= (\hat x + 2^{mn}\hat y)^{\bs(2^n-y)}.
\endalign$$
The isomorphism $\hat\phi$ sends $\hat x + 2^{mn}\hat y$ to
$2^n \hat w + \hat y$, so it sends $(2^{l+n}-2^n+y+1)\bs i$ to
$(a+1)\bs i$; hence, there must exist numbers $v$ and $z$ with $z < 2^n$
such that $(2^{l+n}-2^n+y+1) \bs i = (\hat x + 1)\bs v + z$ and
$(a+1)\bs i = (2^n \hat w + 1)\bs v + z$.  The isomorphism $\phi$
sends $\hat x$ to $2^n \hat w$, so it sends $(\hat x + 1)\bs v =
\hat x^{\bs(2^l-v)}$ to $(2^n \hat w + 1)\bs v = (2^n \hat w)^{\bs(2^l-v)}$,
so there must be a number $u < 2^l$ such that $(\hat x + 1)\bs v =
(x+1) \bs u$ and $(2^n \hat w + 1) \bs v = 2^nu$.  Therefore,
$(a+1)\bs i = (x+1)\bs u + z$ and $(2^{l+n}-2^n+y+1)\bs i = 2^nu + z$,
which is the desired result for $i$.

This completes the case when $x$ and $y$ are both nonzero.  Finally,
we consider the case $y = 0$.  In this case, all that needs to be shown
is that $2^{l+n}-2^n+1$ has maximal period length $2^l$.  But this
follows from Lemma 4.3, which we can apply now since we have proven
the Theorem for $a < L+1$ and nonzero $y$.  This completes the proof
for $a < L+1$, so the induction is complete.
\QED

%% file: critii2.tex
%

\head 5. Applications to algorithms for computing in $A_N$ and $B_N$ \endhead
As noted earlier, computer experimentation in the finite algebras
led to the formulation of Theorem~4.1.  We will now see that this
theorem leads to new algorithms which allow such experimentation to
proceed much farther in the sequence of finite algebras.  Again we
will work with the backward algebras $B_N$; but this clearly gives
algorithms for $A_N$ as well, since the isomorphism between $A_N$ and
$B_N$ is just a trivial two's complement.

Since $B_N$ is defined by recursive formulas, one can easily write
down the obvious recursive algorithm for computing $x\bs y$ in $B_N$.
However, such an algorithm would take a very long time even for small
values of $N$, since very many intermediate products $x'\bs y'$ are
needed, and most of them are needed more than once.  Alternatively,
one could precompute the entire multiplication table for $B_N$ once
relatively quickly and store it, but this requires too much storage space
if $N$ is much more than $10$.

In order to compute $x\bs y$ for $x,y \in B_N$, where $x \ne 0$ (the
case $x=0$ is trivial), it suffices to know the entire period sequence
of $x$; if this period has length $2^l$, then $x\bs y$ is just element
number $(y \bmod 2^l)$ in the period.

Since the period is strictly increasing, there is a threshold $t \le
2^l$ such that, for $0 \le y < 2^l$, we have $x \bs y \ge 2^{N-1}$ if
and only if $y \ge t$.  Hence, if we know $l$ and $t$, we can
immediately read off the $2^{N-1}$\snug's bit of $x \bs y$ for any $y$;
it is $1$ if $y \bmod 2^l \ge t$, $0$ otherwise.  Note that $t > 0$,
since $x \bs 0 = 0$.  If the $2^{N-1}$\snug's bit of $x-1$ is $0$, then the
$2^{N-1}$\snug's bit of $x\bs y$ must be $0$ for any $y$, so $t = 2^l$; if
the $2^{N-1}$\snug's bit of $x \bs y$ is $1$, then $0 < t < 2^l$.

To get the $2^m$\snug's bit of $x\bs y$ for $m \le N-1$, we just have to reduce
everything modulo $2^{m+1}$ and work as above in $B_{m+1}$.  Hence, the
entire period of $x$ can be described by numbers $l(m)$ and $t(m)$ for
$0 \le m < n$; the $2^m$\snug's bit of $x \bs y$ is $1$ iff $y \bmod 2^{l(m)}
\ge t(m)$.

If the $2^m$\snug's bit of $x-1$ is $0$, then one does not need to keep
track of $l(m)$ and $t(m)$, because $l(m)$ is $l(m-1)$ ($0$ if $m = 0$)
and $t(m)$ is $2^{l(m)}$.  Also, if $x'-1 = (x-1) \bmod 2^{n'}$, then
the values $l(m)$ and $t(m)$ for $x'$ are the same as those for $x$ for
all $m < n'$.  This means that, in order to compute $x \bs y$ for all
$x$ and $y$ in $B_N$, it suffices to have a pair of numbers $\bar l(x)$
and $\bar t(x)$ for each $x\ge 2$ in $B_N$, namely the numbers $l(m)$
and $t(m)$ for $x$, where $m$ is the position of the topmost $1$ bit in
$x-1$.

The values of $t(m)$ and $\bar t(x)$ are not arbitrary, as we will show
from the following lemmas.

\proclaim{Lemma 5.1} If $x$ has at most two\/ $1$\snug's in its binary
expansion, then $x+1$ has maximal period length in $\Binf$
(\/$1$, $2$, or\/ $4$
if $x$ has zero, one, or two\/ $1$\snug's, respectively). \endproclaim

\demo{Proof} The case of zero $1$ bits is trivial.  If $x$ has exactly one
$1$ bit, then the period of $x+1$ has length at most $2$ but includes
both $0$ and $x$, so the period length is exactly $2$.  If $x$ has exactly
two $1$ bits, suppose $x = 2^M + 2^m$ with $M>m$; then $2^M+1$ has period
length $2$ by the above, so $x+1$ has period length $4$ by Corollary 4.2
(with $n = 0$). \QED

\proclaim{Corollary 5.2} If $x$ has at least two\/ $1$\snug's in its binary
expansion, then $x+1$ has period length at least\/ $4$. \endproclaim

\demo{Proof} Apply Lemma 5.1 to $x \bmod 2^{m+1}$, where $m$ is the
position of the second-to-last $1$ in the binary expansion of $x$. \QED

\proclaim{Lemma 5.3} For any $x$ and $n$, $x \bs 2^n$ is\/ $0$ or
a power of\/ $2$. \endproclaim

\demo{Proof} We have $(x \bs 2^n)\bs(x \bs 2^n) = x \bs (2^n \bs 2^n) =
x \bs 0 = 0$, so, by Corollary 5.2, $x \bs 2^n$ has at most one $1$ bit. \QED

\proclaim{Lemma 5.4} If $l(m)$ and $t(m)$ are as defined above for $x$,
and $m>0$,
then either the period of~$x$ doubles at $m$ (i.e., $l(m) = l(m-1)+1$)
and $t(m) = 2^{l(m)-1}$, or the period of $x$ does not double at~$m$
(i.e., $l(m) = l(m-1)$) and $t(m) > 2^{l(m)-1}$. \endproclaim

\demo{Proof} The first case is clear, since the period of $(x \bmod 2^{m+1})+1$
is the period of $(x \bmod 2^m)+1$, followed by that same period with
$2^m$ added to each number.  For the second case, it will suffice to
show that, if $l(m) = l(m-1) > 0$, then $(x \bmod 2^{m+1}) \bs 2^{l(m)-1}
< 2^m$.  But $(x \bmod 2^{m+1}) \bs 2^{l(m)-1}$ must be $0$ or a power of $2$
by Lemma 5.3, and it must be nonzero modulo $2^m$ because
$(x \bmod 2^m)\bs 2^{l(m)-1} \ne 0$, so it must be less than $2^m$. \QED

So, if the $2^m$\snug's bit of $x-1$ is $1$, then $2^{l(m)-1} \le t(m)
< 2^{l(m)}$, so $l(m)$ is the least $r$ such that $2^r>t(m)$.

Similar statements can be made if one is working in $A_N$ rather than
$B_N$.  Again periods are strictly increasing and of length a power of
$2$, so a complete period can be specified by numbers $l(m)$ and $t(m)$
for $m < N$.  However, the maximal period length for $x$ is now the
number of $0$\snug's in the expansion of $x$ as an $N$-bit binary
number, rather than the number of $1$\snug's in $x-1$.  If there are at
most two such $0$\snug's, then $x$ has maximal period length; if there
are at least two, then the period length of $x$ is at least $4$.  Lemma
5.3 now states that the midpoint of a period must be of the form
$2^N-2^m$ for some $m$, and Lemma 5.4 states that the threshold $t(m)$
must satisfy $1 \le t(m) \le 2^{l(m)-1}$.  Note that now $t(m) =
2^{l(m)-1}$ if the period has just doubled, but we can also have $t(m)
= 2^{l(m)-1}$ even if the period has not doubled.

To compute in $B_N$, it suffices to know the values ($\bar l(x)$ and)
$\bar t(x)$ for all $x \ge 2$ in $B_N$.  One can precompute these values
recursively; knowing $l(y)$ and $t(y)$ for $y<x$ allows one to do
repeated multiplication by $x-1$ in order to generate the entire period
of $x$, from which $l(x)$ and $t(x)$ can be determined.  Storing these
precomputed values in a file is reasonable for $N$ up to $20$ or slightly
more, so the above gives a reasonable algorithm for computing in $B_N$
for such $N$.

An algorithm that will handle larger algebras is based on the following
lemma.

\proclaim{Lemma 5.5} Suppose\/ $2^N-2^s+1$ has maximal period length
$2^{N-s}$ in $B_N$.  Then, for any $x < 2^{N-s}$, if the period length
of $x+1$ is\/ $2^l$, then the period length of\/ $2^s(x+1)$ is\/ $2^{l+s}$
and\/ $2^s(x+1)\bs(2^si+j) = 2^s((x+1)\bs i) + j$ for all $i<2^l$ and
$j < 2^s$. \endproclaim

\demo{Proof}
Since $2^N-2^s+1$ has maximal period length, $\subalg{2^N-2^s}$ is
isomorphic to $\subalg{2^{N-s}-1}$, and the isomorphism maps $2^sw$ to
$w$ for all $w$, so $2^sw\bs 2^sw' = 2^s(w\bs w')$ for all
$w,w'<2^{N-s}$.

In particular, for any $w < 2^{N-s}$, we have $(2^sw+2^s) \bs (2^N-1) =
2^sw+2^s-1$ and $(2^sw+2^s) \bs (2^N-2^s) = 2^sw$. Since the period of
$2^sw+2^s$ is strictly increasing, and its length is at least $2^s$ (as
one can see by reducing modulo $2^s$), we must have $(2^sw+2^s) \bs
(2^N-z) = 2^sw+2^s-z$ for all $z$, $1 \le z \le 2^s$.  In particular,
this gives $(2^sw+j+1)\bs(2^sw+2^s-1) = 2^sw+j$ for $j < 2^s-1$.

If $x < 2^{N-s}$ and $x'$ is in $\subalg x$, then
$2^sx'$ is in $\subalg{2^sx}$.  We have $2^sx'\bs
(2^sx+2^s-1) = 2^s\hat x + 2^s-1$ for some $\hat x$, and the method of
Lemma 4.8 lets us show that $\hat x = x'\bs x$:  $$2^s\hat x = (2^s\hat
x + 2^s-1)^{\bs(2^s)} = (2^sx' \bs (2^sx+2^s-1)) ^{\bs(2^s)} = 2^sx'
\bs (2^sx+2^s-1)^{\bs(2^s)} = 2^sx' \bs 2^sx = 2^s(x'\bs x).$$  Also,
since $x'$ is in $\subalg x$, its binary
representation has $1$\snug's only where that of $x$ has $1$\snug's, so
$2^sx'+2^s-1$ ends in no more $1$ bits than $2^sx+2^s-1$ does, so
$2^sx+2^s$ is divisible by at least as large a power of~$2$ as
$2^sx'+2^s$ is.  This means that, for any $j < 2^s-1$, since
$(2^sx'+j+1) \bs (2^sx'+2^s-1) = 2^sx'+j$, $2^sx'+2^s$ is divisible by
the period length of $2^sx'+j+1$, so $2^sx+2^s$ is divisible by this
period length, so $(2^sx'+j+1) \bs (2^sx+2^s-1) = 2^sx'+j$.

Given the above facts, an easy downward induction on $a \le 2^N$ shows
that, if $a = 2^si+j$ with $j < 2^s$, then
$2^s(x+1)\bs(2^si+j) = 2^s((x+1)\bs i) + j$.
\QED

Now suppose that we want to compute the period of $a+1$ in $B_N$, and we
have a positive $s$ such that $2^N-2^s+1$ has period length $2^{N-s}$;
suppose also that we already know the period of $b+1$ for all $b$ such that
$2^N-2^s \le b < 2^N$.
Write $a$ in the form $2^sx+y$ where $y < 2^s$.  
If $x = 0$, then we can get the period of $a+1 = y+1$ from the period
of $2^N-2^s+y+1$ by reducing modulo $2^s$ (i.e., just take the values
$l(m)$ and $t(m)$ for $m < s$), so assume $x > 0$.

Let $a_1 = 2^sx+2^s-1$, and
let $2^l$ be the period length of $x$.  Then, by Lemma 5.5, the period
length of $a_1+1$ is $2^{l+s}$, and there is an isomorphism $\phi_1$ between
$\subalg{2^{l+s}}$ and $\subalg{a_1}$ defined by the formula
$\phi_1(2^si+j)=2^s((x+1)\bs i) + j$ for $i < 2^l$ and $j < 2^s$.
Knowing such an isomorphism is equivalent to knowing the period
of $a_1+1$, since $(a_1+1)\bs b = \phi_1(b \bmod 2^{l+s})$.

Let $a_2 = 2^{l+s}-2^s+y$; then $\phi_1(a_2) = a$, so $a_2+1$ and $a+1$
have the same period length.  If $2^{l'}$ is this common period length,
then there is an isomorphism $\phi_2$ from $\subalg{2^{l'}-1}$ to
$\subalg{a_2}$.  Hence, $\phi_1\circ\phi_2$ gives an isomorphism
between $\subalg{2^{l'}-1}$ and $\subalg{\phi_1(a_2)} = \subalg a$, so
the isomorphisms $\phi_1$ and $\phi_2$ (i.e., the periods of $a_1+1$
and $a_2+1$) determine the period of $a+1$.
Another way of looking at
this is that $a = (a_1+1)\bs a_2$, so $a + 1 = (a_1+1)\bs((a_2+1)-1) +
1$, and the right hand side of this equation is just the formula for
$(a_1+1)\circ(a_2+1)$ in the backward version of the algebra $P_N$.

We can determine $\phi_2$ from the (known) period of $2^N-2^s+y+1$ by
just keeping those values $l(m)$ and $t(m)$ for $m < l+s$.  To determine
$\phi_1$ it suffices to determine the period of $x+1$, and this can be
done by a recursive call to this entire procedure.  (Since $x \le a/2^s$,
each level of recursion reduces the argument by a factor of $2^s$, so
not many levels of recursion are needed.)  Finally, given $\phi_1$ and
$\phi_2$, one can compose them to get the period of $a$.

Actually, getting from the $l_i(m),t_i(m)$ forms of $\phi_i$
($i = 1,2$) to the $l(m),t(m)$ form of $\phi_1 \circ \phi_2$ does not
seem to be trivial.  The $l$ part is easy: $l(m) = l_2(l_1(m))$.  But
computing $t(m)$ from $l_1,l_2,t_1,t_2$ seems to require some tricky
Boolean and arithmetical manipulations; no general formula seems apparent.
There are some easy special cases: if $t_1(m) = 2^{l_1(m)}$, then
$t(m) = 2^{l(m)}$; if $t_1(m) = 2^{l_1(m)-1}$ (the period doubling case
for $\phi_1$), then $t(m) = t_2(l_1(m)-1)$; if $l_1(m) = l_1(m-1)$ and
$t_1(m) = t_1(m-1)$, then $t(m) = t(m-1)$.
It turns out that most $m$\snug's fall under these cases;  for those
which do not, the easiest approach seems to be to do a binary
search to find the largest $u < 2^{l(m)}$ such that $\phi_1(\phi_2(u))
< 2^m$, and let $t(m) = u+1$.

Therefore, we have an algorithm for computing in $B_N$ which only requires
storing period data for the last
$2^s$ members of $B_N$, if $2^N-2^s+1$ has maximal
period length.  Lemma 4.3 gives many cases in which this holds:
if $s$ is a power of $2$, then any value of $N$ up to $3s$ is permissible.
In particular, one can work in $B_{48}$ if one has stored period data for
the last $2^{16}$ members of $B_{48}$.

This leaves the problem of generating such period data to be stored and
used.  Fortunately, there is an efficient algorithm for producing such
data by recursion on $N$.  The case $N=0$ is trivial.  Now suppose we
want the data for $N$ and $s$, and we have such data for $N-1$ and $s$.
The data for $N-1$ and $s$ suffice to do operations $(a+1)\bs b$ in
$B_N$ whenever $a < 2^N-2^s$, because the value of $a_2$ arising in
the algorithm above is less than $2^l$ where $l$ is at most the number
of $1$ bits in the binary expansion of $a_1$, and $a<2^N-2^s$ implies
$a_1 < 2^N-1$.  Now we can simply proceed to compute the entire
period of $2^N-2^s+y+1$ for $y = 0,1,\dots ,2^s-1$ successively, by
iterated multiplication by $2^N-2^s+y$.  (The cases $y=0$ and
$y=2^N-1$ should be handled separately, because they give very long
periods which are already known.)

The algorithms described here have been implemented in C++, and work
well for $N=48$ (using a $2.7$-megabyte period data file).  Increasing
$N$ beyond $48$ would require increasing $s$ from $16$ to $32$, and hence
would require a period data file many gigabytes long, so we seem to
have reached the limits of this algorithm.

\head 6. Applications to critical point computations \endhead

First, we translate some of the results about finite algebras into
embedding terms.  In this section, letters $e$ through $k$ will denote
members of $\jalg$ or $\freeone$,
while other letters will denote nonnegative integers.  Also,
we use the notation $k+n$ to denote $k$ applied to $j$ $n$ times:
$k+0 = k$, $k+n+1 = (k+n)(j)$.  (So $\jsub n + m = \jsub{n+m}$.
This matches the notation from $A_N$,
but is the reverse of that from $B_N$.)

\proclaim{Lemma 6.1} Suppose $n$ is a power of\/ $2$ and\/ $\crit k
\ge \gamma_{(m+1)n}$.  Then,
for any $y$ such that\/ $1 \le y \le 2^n$, we have
$k+2^{mn}y \Lequiv(k+2^{(m+1)n}-1)(\gamma_n)/ (k+2^{mn}-1)(\jsub y)$;
equivalently,
$k+2^{mn}y-1 \Lequiv(k+2^{(m+1)n}-1)(\gamma_n)/ (k+2^{mn}-1)\circ \jsub{y-1}$.
\endproclaim

\demo{Proof} Let $N$ be such that $\gamma_N = (k+2^{(m+1)n}-1)(\gamma_n)$,
and apply Corollary 4.2 with $x$ set to the member of $B_N$ corresponding
to $k+2^{(m+1)n}$. To get the second form, use the fact that
$e *_N 1 = e' *_N 1$ implies $e = e'$ in $A_N$ and hence in $P_N$. \QED

\proclaim{Lemma 6.2} If\/ $2^n \mid s$, then $\jsub{2^s-1}(\gamma_m) =
\gamma_{s+m}$ for all $m < 2^{n+1}$.
\endproclaim

\demo{Proof} Apply Lemma 4.3 with $s' = s+2^{n+1}$. \QED

\proclaim{Lemma 6.3} Suppose $n$ is a power of\/ $2$, $s$ is a multiple of
$n$, and\/ $\crit k \ge \gamma_s$.  Then,
for any $y$ such that\/ $1 \le y \le 2^s$, we have
$k+y \Lequiv(k+2^s-1)(\gamma_n)/ k(\jsub y)$; equivalently,
$k+y-1 \Lequiv(k+2^s-1)(\gamma_n)/ k \circ \jsub{y-1}$. \endproclaim

\demo{Proof} Let $N$ be such that $\gamma_N = (k+2^s-1)(\gamma_n)$,
and apply Lemma 4.4 with $x$ set to the member of $B_N$ corresponding
to $k+2^s$. \QED

\proclaim{Corollary 6.4} If $n$, $s$, and $k$ are as in Lemma~6.3,
then $k(\gamma_{n+s}) \ge (k+2^s-1)(\gamma_n)$. \endproclaim

\demo{Proof} Let $\delta = (k+2^s-1)(\gamma_n)$; then, by
Lemma~6.3, $k+2^s-1 \Lequiv\delta/ k \circ \jsub{2^s-1}$.  Since
$(k+2^s-1)(\gamma_n) \ge \delta$, we must have $(k \circ
\jsub{2^s-1})(\gamma_n) \ge \delta$.  Lemma~6.2 gives
$\jsub{2^s-1}(\gamma_n) = \gamma_{n+s}$, so we have
$k(\gamma_{n+s}) \ge \delta$, as desired. \QED

\proclaim{Lemma 6.5} If $n$ is a power of\/ $2$, $\crit k \ge \gamma_{(m+1)n}$,
and\/ $2 \le z \le 2^n$, then $(k+2^{mn}z-1)(\gamma_n) \ge (k+2^{(m+1)n}-1)
(\gamma_n)$, and hence $(k+2^{mn}z-1)(\gamma_n) \ge \gamma_{(m+2)n}$.
\endproclaim

\demo{Proof} Let $N$ be such that $\gamma_N = (k+2^{(m+1)n}-1)(\gamma_n)$,
and apply Lemma 4.6 with $x$ set to the member of $B_N$ corresponding
to $k+2^{(m+1)n}$ to get the first part. For the second part, note that
$(k+2^{(m+1)n}-1)(\gamma_0) \ge \gamma_{(m+1)n}$ and hence
$(k+2^{(m+1)n}-1)(\gamma_n) \ge \gamma_{(m+2)n}$. \QED

\proclaim{Proposition 6.6} If $k \bmod 2^N$ has at least one\/ $0$ in
its $N$-bit binary representation, then $k(\gamma_0) = \gamma_n$ where
$n$ is the position of the rightmost such\/ $0$ (i.e., $n$ is least
such that $k \bmod 2^{n+1} < 2^n$).  If $k \bmod 2^N$ has at least
two\/ $0$\snug's in its $N$-bit binary representation, then $k(\gamma_1) =
\gamma_{n'}$ where $n'$ is the position of the second-rightmost such\/ $0$.
\endproclaim

\demo{Proof} Apply Lemma 5.1 with $x+1$ equal to the element of
$B_N$ corresponding to $k$. \QED

Next, we give an extension of Lemma 6.1 which is a variant of Lemma~2.8.

\proclaim{Lemma 6.7} Suppose $n$ is a power of\/ $2$ and\/ $\crit k
\ge \gamma_{(m+1)n}$.  Then,
for any $y$ such that\/ $1 \le y \le 2^n$, we have
$k+2^{mn}y \Lequiv\theta_y/ (k+2^{mn}-1)(\jsub y)$, where
$\theta_y = \min_{1 < z \le y} (k+2^{mn}z-1)(\gamma_n)$.
\endproclaim

\demo{Proof} Induct on $y$. The case $y=1$ is trivial, so suppose
$y > 1$.  We then have
$$(k+2^{mn}-1)(\jsub y) = (k+2^{mn}-1)(\jsub{y-1})(k+2^{mn})
\Lequiv\theta_{y-1}/ (k+2^{mn}(y-1))(k+2^{mn}).$$
We have $k+2^{mn} \Lequiv\gamma_{(m+1)n}/ \jsub{2^{mn}}$, so
$$(k+2^{mn}(y-1))(k+2^{mn}) \Lequiv (k+2^{mn}(y-1))(\gamma_{(m+1)n})/
(k+2^{mn}(y-1))(\jsub{2^{mn}}).$$
Apply Lemma 6.3 with $s$, $k$, and $y$ replaced by $mn$,
$k+2^{mn}(y-1)$, and $2^{mn}$ respectively to get
$k+2^{mn}y \Lequiv(k+2^{mn}y-1)(\gamma_n)/ (k+2^{mn}(y-1))(\jsub{2^{mn}})$.
Also by Lemma 6.3, we have $$k+2^{mn}y-1 \Lequiv(k+2^{mn}y-1)(\gamma_n)/
(k+2^{mn}(y-1)) \circ \jsub{2^{mn}-1};$$ since Lemma 6.2 gives
$\jsub{2^{mn}-1}(\gamma_n) = \gamma_{(m+1)n}$, we get
$$(k+2^{mn}(y-1))(\gamma_{(m+1)n}) = ((k+2^{mn}(y-1)) \circ \jsub{2^{mn}-1})
(\gamma_n) \ge (k+2^{mn}y-1)(\gamma_n).$$
Therefore, $\theta_y$ is the minimum of $\theta_{y-1}$,
$(k+2^{mn}(y-1))(\gamma_{(m+1)n})$, and $(k+2^{mn}y-1)(\gamma_n)$,
so $k+2^{mn}y \Lequiv\theta_y/ (k+2^{mn}-1)(\jsub y)$, as desired.
\QED

This works even if $k$ is the identity (which is thought of as having
critical point greater than any~$\gamma_N$), so that
$k+r$ is just $\jsub r$.

\head 7. Column 2 \endhead
Column~2, in both the regular construction and the irregular construction,
is a list
of embeddings $b_0,b_1,\dotsc$ such that \threeseq
b_n/\zeta'/\beta_n/\beta_{n+1}/, where we also have available
an embedding $z$ such that \threeseq
z/\zeta/\zeta'/\beta_0/. Given a list of this type, one can generate
another list of the type suitable for column~1, which results in the
production of $2^{2^n}$ critical points below~$\beta_n$
\cite{\Dougherty,~Lemma~3}.  Hence, if such a list
is to be a column of the `optimal' construction, we must have
$\beta_n = \gamma_{2^{2^n}}$ (and $\zeta = \gamma_0$, $\zeta' = \gamma_1$).
We will see in this section that a suitable sequence of embeddings
does exist.  The proof will be similar to that of Theorem~3.1, but
will require the results of the previous section, as well as the
following lemma.

\proclaim{Lemma 7.1} Given $k \in \jalg$, let $f$ be the numerical
mapping associated with $k$, so that $k(\gamma_n) = \gamma_{f(n)}$.
Then $f(n) \ge f(0) + (f(1)-f(0))n$ for all $n$.  Furthermore,
if we have $f(n) = f(0) + n(f(1)-f(0))$ for a particular $n$, then
$f(n') = f(0) + (f(1)-f(0))n'$ for all $n' \le n$, and
$k \bmod 2^{f(n)+1} < 2^{f(1)}$. \endproclaim

\demo{Proof} As mentioned at the beginning of section~3, for
any $m>0$, \twoseq\jsub{2^m}/\gamma_m/\gamma_{m+1}/.  But
\twoseq j/\gamma_0/\gamma_1/, and $\jsub{2^m} = \jsub{2^m-1}(j)$,
so we must have $\jsub{2^m-1}(\gamma_0) = \gamma_m$ and
$\jsub{2^m-1}(\gamma_1) = \gamma_{m+1}$.  Now apply $k$ to this to get
$k\jsub{2^m-1}(\gamma_{f(0)}) = \gamma_{f(m)}$ and
$k\jsub{2^m-1}(\gamma_{f(1)}) = \gamma_{f(m+1)}$.
Hence, $k\jsub{2^m-1}$ maps the critical points between $\gamma_{f(0)}$
and~$\gamma_{f(1)}$ to critical points between $\gamma_{f(m)}$
and $\gamma_{f(m+1)}$, so we must have $f(m+1)-f(m) \ge f(1)-f(0)$.
This is clearly true for $m=0$ as well.
Adding these inequalities for $m=0,1,\dots,n-1$ gives
$f(n)-f(0) \ge (f(1)-f(0))n$, as desired.

If $f(n)-f(0) = (f(1)-f(0))n$, then we must have had
$f(m+1)-f(m) = f(1)-f(0)$ for all $m<n$, so we get
$f(n')-f(0) = (f(1)-f(0))n'$ for all $n'\le n$.
We now show that, for $f(1) \le r \le f(n)$, bit number $r$ from the right
(the $2^r$\snug's bit) of $k \bmod 2^{f(n)+1}$ is $0$;
this will give $k \bmod 2^{f(n)+1} < 2^{f(1)}$.
If $r$~is $f(m)$ for some $m$, then $k(\gamma_m)=\gamma_r$, so
Lemma~2.9 states that bit number $r$ of $k \bmod 2^{f(n)+1}$ must be~$0$.
Now suppose $r$ is not one of the numbers $f(m)$; then $r$~must lie
between $f(m)$ and $f(m+1)$ for some $m$, with $1 \le m < n$.
Since $f(m+1)-f(m) = f(1)-f(0)$, $k\jsub{2^m-1}$ must map the
entire list of critical points between $\gamma_{f(0)}$ and $\gamma_{f(1)}$
to the entire list of critical points between $\gamma_{f(m)}$ and
$\gamma_{f(m+1)}$.  In particular, we must have
$k\jsub{2^m-1}(\gamma_s)=\gamma_r$ for some $s$, so
bit number~$r$ of $k\jsub{2^m-1} \bmod 2^{f(n)+1}$ is~$0$.
Since $\crit{(k\jsub{2^m-1})} = \gamma_{f(0)}$, we have
$k\jsub{2^m-1} \bmod 2^{f(0)+1} = 2^{f(0)}$, so
bit number~$f(0)$ of $k\jsub{2^m-1} \bmod 2^{f(n)+1}$ is~$1$.
This means that subtracting~$1$ from $k\jsub{2^m-1} \bmod 2^{f(n)+1}$
does not require borrowing from bit number~$r$, so
bit number~$r$ of $(k\jsub{2^m-1} \bmod 2^{f(n)+1}) - 1$ is~$0$.
Therefore, by the remarks preceding Lemma~2.9,
bit number~$r$ of $k \bmod 2^{f(n)+1}$ is~$0$, as desired. \QED

We will also need results guaranteeing the existence of embeddings whose critical sequences begin
in particular ways.  One such result is Theorem~3.1(a).  Another is:
For any $n \ge 2$, there exists $k \in \jalg$ such that
\threeseq k/\gamma_0/\gamma_1/\gamma_n/.  This is proved
by induction on $n$.  For the base case, since $\kappa_2 = \gamma_2$
(this follows from $\crit{\jsub 4} = \kappa_2$, which is proved
by an easy computation~\cite{\LaverII, \Dougherty}),
\threeseq j/\gamma_0/\gamma_1/\gamma_2/.  For the induction
step, if \threeseq k/\gamma_0/\gamma_1/\gamma_n/, then
\threeseq \jsub{2^n}k/\gamma_0/\gamma_1/\gamma_{n+1}/.

\proclaim{Theorem 7.2} For any natural number $n$, the following are true:
\roster
\item"\rm(a)" For any $m < n$, there exists $\kk$ in $\jalg$
such that
\threeseq\kk/\gamma_{2^{2^n-2^{m+1}}}/\gamma_{2^{2^n-2^m}}/\gamma_{2^{2^n}}/.
\item"\rm(b)" If $s$ is a multiple of $2^{2^n}$, and $k \bmod
2^{s+2^{2^n}} = 2^sx+2^s-1$ where $1 \le x < 2^{2^{2^n}} - 1$, then
$k(\gamma_{2^{2^n}}) \ge \gamma_{s+2^{2^{n+1}}}$.  Furthermore, if also
$k \bmod 2^{s+2^{2^{n+1}}} > 2^{s+2^{2^n}}$, then $k(\gamma_{2^{2^n}})
> \gamma_{s+2^{2^{n+1}}}$.
\item"\rm(c)"
\threeseq\jsub{2^{2^{2^n}}-2}/\gamma_1/\gamma_{2^{2^n}}/\gamma_{2^{2^{n+1}}}/.
\endroster \endproclaim

\demo{Proof} By simultaneous induction on $n$.  The case $n=0$ of
(a) is vacuously true.  Throughout this proof, let $N = 2^{2^n}$; then
$N^2 = 2^{2^{n+1}}$.

First, suppose (a) is true for $n$.  We prove (b) for $n$ by downward
induction on $x$.

If $x = 2^{N}-2$, then Proposition~6.6 gives $k(\gamma_0) = \gamma_s$
and $k(\gamma_1) \ge \gamma_{s+N}$.  Therefore, Lemma~7.1
implies $k(\gamma_r) \ge \gamma_{s+Nr}$ for all $r$, so
$k(\gamma_{N}) \ge \gamma_{s+N^2}$.  Furthermore, if
we have equality here, then Lemma~7.1 implies $k(\gamma_1) = \gamma_{s+N}$
and $k \bmod 2^{s+N^2} < 2^{s+N}$.

Now suppose $1 \le x < 2^{N}-2$.  Then Proposition~6.6 gives $\gamma_s
\le k(\gamma_0) < \gamma_{s+N}$, and, if $k(\gamma_0) = \gamma_s$, then
$k(\gamma_1) < \gamma_{s+N}$.

\procl{Claim 1} There exist critical points $\alpha$ and $\beta$ and
an embedding $\kk \in \jalg$ such that
$\gamma_s < k(\alpha) < \gamma_{s+N} \le k(\beta)$ and
\threeseq \kk/\alpha/\beta/\gamma_{N}/. \endproclaim

\demo{Proof} We consider several cases.

Case 1: $k(\gamma_1) \ge \gamma_{s+N}$.  Then, since
$x < 2^{N}-2$, we must have
$\gamma_s < k(\gamma_0) < \gamma_{s+N}$.  Let $\alpha = \gamma_0$
and $\beta = \gamma_1$, and obtain $\kk$ as in the paragraph
preceding the statement of Theorem~7.2.

Case 2: $k(\gamma_1) < \gamma_{s+N}$ but
$k(\gamma_{2^{2^n-1}}) \ge \gamma_{s+N}$.  Then there is a
largest $m$ such that $k(\gamma_{2^{2^n-2^m}}) \ge
\gamma_{s+N}$, and this $m$ is less than $n$.
Let $\alpha = \gamma_{2^{2^n-2^{m+1}}}$ and $\beta = \gamma_{2^{2^n-2^m}}$,
and apply (a) for $n$ and $m$ to get $\kk$.

Case 3: $k(\gamma_{2^{2^n-1}}) < \gamma_{s+N}$.  Since
$k \bmod 2^{s+N}$ has at most $N-1$ $0$\snug's among its lowest
$s+N$ bits, we have $k(\gamma_{N-1}) \ge
\gamma_{s+N}$.  Hence, there is a
largest $m$ such that $k(\gamma_{N-2^m}) \ge
\gamma_{s+N}$, and this $m$ is less than $2^n-1$.
Let $\alpha = \gamma_{N-2^{m+1}}$ (this is not $\gamma_0$
because $m < 2^n-1$) and $\beta = \gamma_{N-2^m}$,
and apply Theorem~3.1(a) for $2^n$ and $m$ to get $\kk$.
\QEd

Now $\crit{(k\kk)} = k(\crit \kk) = k(\alpha)$ is strictly between
$\gamma_s$ and $\gamma_{s+N}$, so $k\kk \bmod 2^{s+N} = 2^sx'$
where $1 \le x' \le 2^N-1$ and $x'$ is even (so $x' < 2^N-1$).
By \thetag{2.2}, we have $x' > x$.
Therefore, we can apply the inner
inductive hypothesis to get $(k\kk+2^s-1)(\gamma_{N}) \ge
\gamma_{s+N^2}$.  But Corollary~6.4 gives
$k\kk(\gamma_{s+N}) \ge (k\kk+2^s-1)(\gamma_{N})$, so
$$k(\gamma_{N}) = k(\kk(\beta)) =
k\kk(k(\beta)) \ge k\kk(\gamma_{s+N}) \ge (k\kk+2^s-1)(\gamma_{N}).$$
Therefore, $k(\gamma_{N})
\ge \gamma_{s+N^2}$.  If we also have $k \bmod 2^{s+N^2} >
2^{s+N}$, then \thetag{2.2} gives $k\kk \bmod
2^{s+N^2} > k \bmod 2^{s+N^2} > 2^{s+N}$, so the inner
inductive hypothesis gives $(k\kk+2^s-1)(\gamma_{N}) >
\gamma_{s+N^2}$ and hence $k(\gamma_{N}) >
\gamma_{s+N^2}$.  This completes the proof of (b) for $n$.

Next, we prove (c) for $n$ from (b) for $n$.  Previous results
easily show that \twoseq\jsub{2^N-2}/\gamma_1/\gamma_N/ (Theorem~2.1
shows that $\crit{\jsub{2^N-2}} = \gamma_1$, and
$\jsub{2^N-2}(j(\gamma_0)) = \jsub{2^N-1}(\gamma_0) = \gamma_N$),
so we must show that $\jsub{2^N-2}(\gamma_N)=\gamma_{N^2}$.

Define a numerical function
$h$ as follows: if the number $r$ has binary expansion $\sum_m b_m 2^m$,
where $b_m$ is $0$ or~$1$, then let $h(r) = \sum_m b_m(2^{(m+1)N} - 2^{mN})$.
In other words, $h(r)$ is obtained from $r$ by replacing each $0$ or $1$
bit in the binary expansion of~$r$ with a block of $N$ $0$\snug's or
$1$\snug's.

\procl{Claim 2} For $1 \le r < 2^N-1$, we have $\jsub{h(r+1)}
\Lequiv\gamma_{N^2+1}/ \jsub{h(r)}\jsub{2^N-1}$. \endproclaim

\demo{Proof} Let $m$ be such that $2^m \parallel r+1$ (that is, $m$ is
the location of the rightmost $1$ bit in $r+1$).  Let $y = 2^N-1$
and let $k$ be $\jsub{h(r+1-2^m)}$ (the identity if $r+1=2^m$).
Since $2^m \parallel r+1$, we have $2^{m+1} \mid r+1-2^m$,
so $2^{(m+1)N} \mid h(r+1-2^m)$, so $\crit k \ge \gamma_{(m+1)N}$.
Also, $k+2^{mN}y = \jsub{h(r+1)}$ and $k+2^{mN}-1 = \jsub{h(r)}$.
Now Lemma~6.7 states that $\jsub{h(r+1)}
\Lequiv\theta_y/ \jsub{h(r)}\jsub{N-1}$, where
$\theta_y = \min_{1 < z \le y} (k+2^{mN}z-1)(\gamma_N)$.
But (b) for $n$ gives $(k+2^{mN}z-1)(\gamma_N) > \gamma_{N^2}$ for
all such $z$, so $\theta_y \ge \gamma_{N^2+1}$, and we are done.
\QEd

Apply this for
all $r$ from $2^N-2$ down to $1$, noting that $h(2^N-1) = 2^{N^2}-1$
and $h(1) = 2^N-1$, to get
$$\align \jsub{2^{N^2}-1}
&\Lequiv\hphantom{\gamma_{N^2+1}}/ \jsub{h(2^N-1)} \\
&\Lequiv\gamma_{N^2+1}/ \jsub{h(2^N-2)}\jsub{2^N-1} \\
&\Lequiv\gamma_{N^2+1}/
   \jsub{h(2^N-3)}\jsub{2^N-1}\jsub{2^N-1} \\
&\Lequiv\gamma_{N^2+1}/ \dotsm \\
&\Lequiv\gamma_{N^2+1}/ \jsub{2^N-1}\jsub{2^N-1}\dots\jsub{2^N-1}
   \qquad\text{($2^N-1$ times)} \\
&\Lequiv\hphantom{\gamma_{N^2+1}}/ \jsub{2^N-2}\jsub{2^N-1}.
\endalign$$
Since $\jsub{2^N-2}(\gamma_0) = \gamma_0$, we have
$\jsub{2^N-2}\jsub{2^N-1}(\gamma_0) =
\jsub{2^N-2}\jsub{2^N-1}(\jsub{2^N-2}(\gamma_0)) =
\jsub{2^N-2}(\gamma_N)$.  On the other hand,
$\jsub{2^{N^2}-1}(\gamma_0) = \gamma_{N^2}$.
Since $\jsub{2^{N^2}-1}$ and $\jsub{2^N-2}\jsub{2^N-1}$ agree up to
$\gamma_{N^2+1}$, we must have $\jsub{2^N-2}(\gamma_N) = \gamma_{N^2}$.
Therefore, (c) holds for $n$.

Finally, assuming (a) and (c) hold for $n$, we get (a) for $n+1$ as
follows.  Given (c) for $n$, we can apply Lemma~7.1 to $\jsub{2^{N}-2}$
to get $\jsub{2^{N}-2}(\gamma_m) = \gamma_{Nm}$ for all
$m \le N$.  Now,
to get $a$ for $n+1$ and $m$ where $m = n$, use the embedding
$\jsub{2^{N}-2}$.  For $m < n$, apply (a) for $n$ and $m$ to get
$\kk$; then $\jsub{2^{N}-2}(\kk)$ works for $n+1$ and $m$. \QED

So one may expect column 2 of the ultimate construction
to be a sequence
of embeddings $b_0,b_1,\dotsc$ such that \threeseq b_n/\gamma_1/
\gamma_{2^{2^n}}/\gamma_{2^{2^{n+1}}}/.

\head 8. Column 3? \endhead
Since the first two columns of the regular and irregular constructions
have worked out to be potentially
optimal, it is natural to ask whether the third column does also.
The third column in these constructions is a list
of embeddings $d_0,d_1,\dotsc$ such that \fourseq d_n/\zeta/\zeta'/
\delta_n/\delta_{n+1}/.  Given a list of $n$ embeddings like this,
one can generate a list of $C_2(n)$ embeddings of the column 2 type,
where $C_2(n)$ is defined recursively by $C_2(0) = 0$ and
$C_2(n+1) = C_2(n) + 2^{C_2(n)}$~\cite{\Dougherty,~Lemma~4}.  Hence,
there are at least $2^{2^{C_2(n)}}$ critical points below $\delta_n$.
If such a column is to be optimal, we must have $\delta_n =
\gamma_{2^{2^{C_2(n)}}}$.  Also, note that member $b_{C_2(n)}$ of
column 2 is obtained by applying member $d_n$ of column 3 to the first
member of column 1, which is just $j$; hence, if $b_{C_2(n)}$ `looks
like' $\jsub m$, then $d_n$ should `look like' $\jsub{m-1}$.  So the
third-column situation seems likely to be as follows.

\proclaim{Conjecture 8.1} For all $n$, \fourseq \jsub{2^{f(n)}-3}/
\gamma_0/\gamma_1/\gamma_{f(n)}/\gamma_{f(n+1)}/,
where $f(n) = 2^{2^{C_2(n)}}$.
\endproclaim

For $n=0$, this conjecture states that \fourseq j/
\gamma_0/\gamma_1/\gamma_2/\gamma_4/, which we already know to be true.
We will now see that the conjecture is true for $n=1$ as well;
that is, \fourseq \jsub{13}/
\gamma_0/\gamma_1/\gamma_4/\gamma_{256}/.  Unfortunately, the
methods to be used do not have an obvious generalization to
later members of the column.

We will be using a number of results from section~5 of
Dougherty~\cite{\Dougherty}; for instance, it is shown there that
\threeseq \jsub{13}/\gamma_0/\gamma_1/\gamma_4/.  Also, we have
$\kappa_3 = \gamma_4$, and $\jsub{13}(\kappa_2) = \jsub{14}(\kappa_3)$.
These can be combined with results in the present paper, which
state that $\jsub{14}(\gamma_4) = \gamma_{16}$, to give
$\jsub{13}(\kappa_2) = \gamma_{16}$.

The proofs here will be written in terms of ordinary restrictions
$k \restrict V_\alpha$, rather than Laver restrictions $k \Lrestrict V_\alpha$
and the corresponding equivalence relation $\Lequiv\alpha/$.
This is not one of the methods that automatically works in the finite
algebra context as well because of Proposition~2.7.  However,
this particular proof can be translated in a systematic way into
a proof using Laver restrictions, which does work in the finite
algebra context; for an example of this, see the two proofs of
Lemma~15 in Dougherty~\cite{\Dougherty}.  But the translated proof
is substantially longer and less readable.

If $k(\alpha) = \beta$, and $\beta$ is a limit ordinal, then clearly
$k \restrict V_\alpha \subseteq V_\beta$.  However, if $\beta$ is
inaccessible (as all critical points are) and $\alpha < \beta$,
then the size of $k \restrict V_\alpha$ is less than $\beta$,
so we actually get $k \restrict V_\alpha \in V_\beta$.  For
instance, $j \restrict V_{\kappa_1} \in V_{\kappa_2}$ and
$j \restrict V_{\kappa_2} \in V_{\kappa_3}$.

Straightforward computation~\cite{\Dougherty,~\S5} gives the
following inequalities (which are not in the strong\-est possible form):
$$\alignedat 3
j(\kappa_2) &\ge \kappa_3 &\qquad\qquad \jsub5(\kappa_1) &\ge \kappa_2
   &\qquad\qquad \jsub9(\kappa_1) &\ge \kappa_2 \\
\jsub2(\kappa_2) &\ge \kappa_3 &\qquad\qquad \jsub6(\kappa_1) &\ge \kappa_2
   &\qquad\qquad \jsub{10}(\kappa_1) &\ge \kappa_2 \\
\jsub3(\kappa_2) &\ge \kappa_3 &\qquad\qquad \jsub7(\kappa_1) &\ge \kappa_3
   &\qquad\qquad \jsub{11}(\kappa_1) &\ge \kappa_3 \\
\jsub4(\kappa_2) &\ge \kappa_2 &\qquad\qquad \jsub8(\kappa_2) &\ge \kappa_2
   &\qquad\qquad \jsub{12}(\kappa_2) &\ge \kappa_3
\endalignedat\tag{8.1}$$
We also have $\jsub{13}(\kappa_2)=\gamma_{16}$ (as noted above) and
$\jsub{14}(\kappa_2)=\gamma_{8}$ (from Theorem~7.2 and Lemma 7.1).
We can use these facts to expand out part of $\jsub{14}$,
as follows.  Start with $\jsub{14}\restrict V_{\gamma_{16}} =
\jsub{13}(j\restrict V_{\kappa_2}) = \jsub{12}j(j\restrict V_{\kappa_2})$.
Since $\jsub{12}(\kappa_2) \ge \kappa_3$ and
$j\restrict V_{\kappa_2} \in V_{\kappa_3}$, we actually have
$\jsub{13}(j\restrict V_{\kappa_2}) = \jsub{12}(j\restrict V_{\kappa_2})
(j\restrict V_{\kappa_2})$.  Similarly, since $\jsub{11}(\kappa_1)
\ge \kappa_3$, we get
$\jsub{12}(j\restrict V_{\kappa_2}) = \jsub{11}(j\restrict V_{\kappa_1})
(j\restrict V_{\kappa_2})$.  Continuing this leads to the following
expansion:
$$\multline
\jsub{14}\restrict V_{\gamma_{16}} = j(j\restrict V_{\kappa_2})
(j\restrict V_{\kappa_2}) (j\restrict V_{\kappa_2}) (j\restrict
V_{\kappa_2}) (j\restrict V_{\kappa_1}) (j\restrict V_{\kappa_1})
(j\restrict V_{\kappa_1})\\
(j\restrict V_{\kappa_2}) (j\restrict
V_{\kappa_1}) (j\restrict V_{\kappa_1}) (j\restrict V_{\kappa_1})
(j\restrict V_{\kappa_2}) (j\restrict V_{\kappa_2}).
\endmultline\tag 8.2$$

Lemma~6.3 implies the following: If $\crit k \ge \gamma_{4m}$, then
$k+2^{4m}-1 \Lequiv(k+2^{4m}-1)(\gamma_4)/ k \circ \jsub{2^{4m}-1}$.
Since $j \restrict V_{\kappa_2} \in V_{\gamma_4}$, we can conclude
that $$k(\jsub{2^{4m}-1}(j\restrict V_{\kappa_2})) =
(k+2^{4m}-1)(j \restrict V_{\kappa_2}) = (k+2^{4m}) \restrict
V_{(k+2^{4m}-1)(\kappa_2)}.\tag 8.3$$
If we make the stronger assumption that $\crit k \ge \gamma_{4m+4}$,
then we get a stronger conclusion: Lemma 6.2 gives
$\jsub{2^{4m}-1}(\gamma_4) = \gamma_{4m+4}$, so
$\jsub{2^{4m}-1}(j\restrict V_{\kappa_2}) \in V_{\gamma_{4m+4}}$,
so $$(k+2^{4m}-1)(j \restrict V_{\kappa_2}) =
k(\jsub{2^{4m}-1}(j\restrict V_{\kappa_2})) =
\jsub{2^{4m}-1}(j\restrict V_{\kappa_2}) =
\jsub{2^{4m}}\restrict V_{\gamma_{4m+2}}. \tag 8.4$$
These equalities can then be restricted to smaller domains, to
get $k(\jsub{2^{4m}-1}(j\restrict V_{\kappa_1})) =
(k+2^{4m}) \restrict
V_{(k+2^{4m}-1)(\kappa_1)}$ from~\thetag{8.3},
and similarly for~\thetag{8.4}.

\proclaim{Lemma 8.2} If $\crit k \ge 4m+4$, then
$(k+2^{4m}-1)(\jsub{14}\restrict V_{\gamma_{16}}) =
(k+14\cdot2^{4m})\restrict V_{(k+2^{4m}-1)(\gamma_{16})}$ and
$(k+2^{4m}-1)(\jsub{15}\restrict V_{\gamma_8}) =
(k+15\cdot2^{4m})\restrict V_{(k+2^{4m}-1)(\gamma_8)}$. \endproclaim

\demo{Proof} Plug in \thetag{8.2}, use the distributive law,
and then apply \thetag{8.4} and \thetag{8.3} repeatedly:
$$\align
&(k+2^{4m}-1)(\jsub{14}\restrict V_{\gamma_{16}}) \\
&\qquad= (k+2^{4m}-1)(j(j\restrict V_{\kappa_2})
(j\restrict V_{\kappa_2}) \dotsm (j\restrict V_{\kappa_2})) \\
&\qquad= (k+2^{4m})((k+2^{4m}-1)(j\restrict V_{\kappa_2}))
((k+2^{4m}-1)(j\restrict V_{\kappa_2})) \dotsm
((k+2^{4m}-1)(j\restrict V_{\kappa_2})) \\
&\qquad= (k+2^{4m})(\jsub{2^{4m}-1}(j\restrict V_{\kappa_2}))
(\jsub{2^{4m}-1}(j\restrict V_{\kappa_2})) \dotsm
(\jsub{2^{4m}-1}(j\restrict V_{\kappa_2})) \\
&\qquad= (k+2\cdot2^{4m})
(\jsub{2^{4m}-1}(j\restrict V_{\kappa_2})) \dotsm
(\jsub{2^{4m}-1}(j\restrict V_{\kappa_2})) \\
&\qquad= (k+3\cdot2^{4m}) \dotsm
(\jsub{2^{4m}-1}(j\restrict V_{\kappa_2})) \\
&\qquad= \dotsm \\
&\qquad= (k+13\cdot2^{4m}) 
(\jsub{2^{4m}-1}(j\restrict V_{\kappa_2})) \\
&\qquad= (k+14\cdot2^{4m}) \restrict V_{(k+13\cdot2^{4m})(\kappa_{4m+2})}
\endalign$$
But clearly two functions which are equal must have the same domain,
so $(k+2^{4m}-1)(V_{\gamma_{16}}) = V_{(k+13\cdot2^{4m})(\kappa_{4m+2})}$
and we get the desired result.  The argument for $\jsub{15}$
is the same with one more $(j\restrict V_{\kappa_2})$. \QED

We now use the function $h$ defined in the proof of Theorem~7.2,
letting $N=4$: $h(r)$ is obtained from $r$ by quadruplicating
every bit in the binary expansion of $r$.  So
$h(1) = 15$, $h(2) = 240$, $h(14) = 2^{16}-16$, and so on.

We can now show that, for any $r$,
$\jsub{h(r)}(\jsub{15}\restrict V_{\gamma_8}) =
\jsub{h(r+1)} \restrict V_{\jsub{h(r)}(\gamma_8)}$.  To see this,
let $m$ be such that $2^m \parallel r+1$, and apply Lemma~8.2
with $k = \jsub{h(r+1-2^m)}$.  Again, if we want, we can simply restrict
both sides of this equality further to get
$\jsub{h(r)}(\jsub{15}\restrict V_{\gamma_4}) =
\jsub{h(r+1)} \restrict V_{\jsub{h(r)}(\gamma_4)}$.

Now we can use \thetag{8.2} again, to get
$$\align
\jsub{13}(jj)\restrict V_{\jsub{13}(\kappa_3)}
&= \jsub{13}(j(j\restrict V_{\kappa_2})) \\
&= \jsub{14}(\jsub{14}\restrict V_{\gamma_{16}}) \\
&= \jsub{14}(j(j\restrict V_{\kappa_2})
(j\restrict V_{\kappa_2}) \dotsm (j\restrict V_{\kappa_2})) \\
&= \jsub{h(1)}(\jsub{15}\restrict V_{\gamma_8}) (\jsub{15}\restrict
V_{\gamma_8}) \dotsm (\jsub{15}\restrict V_{\gamma_8}) \\
&= \jsub{h(2)}(\jsub{15}\restrict
V_{\gamma_8}) \dotsm (\jsub{15}\restrict V_{\gamma_8}) \\
&= \jsub{h(3)} \dotsm (\jsub{15}\restrict V_{\gamma_8}) \\
&= \dotsm \\
&= \jsub{h(13)} (\jsub{15}\restrict V_{\gamma_8}) \\
&= \jsub{h(14)} \restrict V_{\jsub{h(13)}(\gamma_8)}.
\endalign$$
Hence, by another application of Lemma~8.2,
$$\align
\jsub{13}(\jsub3)\restrict V_{\jsub{13}(\kappa_3)}
&= \jsub{13}(j(j\restrict V_{\kappa_2})(j\restrict V_{\kappa_2})) \\
&= \jsub{13}(j(j\restrict V_{\kappa_2}))(\jsub{14}\restrict V_{\gamma_{16}}) \\
&= \jsub{h(14)}(\jsub{14}\restrict V_{\gamma_{16}}) \\
&= \jsub{2^{16}-2}\restrict V_{\jsub{h(14)}(\gamma_{16})}.
\endalign$$
We can now apply this to $\gamma_{16}$ to get
$$\jsub{2^{16}-2}(\gamma_{16}) = \jsub{13}(\jsub3)(\gamma_{16})
= \jsub{13}(\jsub3)(\jsub{13}(\kappa_2)) = \jsub{13}(\kappa_3).$$
But Theorem~7.2 gives $\jsub{2^{16}-2}(\gamma_{16}) = \gamma_{256}$,
so we have $\jsub{13}(\kappa_3) = \gamma_{256}$, as desired.

\medskip
We can now use this evaluation of $\jsub{13}(\kappa_3)$ to compute
some other critical points.  For instance, since $\jsub7(\kappa_2)$,
$\jsub{12}(\kappa_3)$, and $\jsub{13}(\kappa_3)$ are equal to
each other~\cite{\Dougherty, \S5}, they are all equal to $\gamma_{256}$.
A more interesting fact is as follows:

\proclaim{Proposition 8.3} The\/ $256$ critical points below $\kappa_4$
produced by the regular construction\/~{\rm\cite{\Dougherty,~\S2}} are
the first\/ $256$ critical points $\gamma_0,\gamma_1,\dots,\gamma_{255}$.
\endproclaim

\demo{Proof} Suppose we modify the regular construction by
replacing the first entry in column~4, which is normally $(jj)(jj) =
j\jsub2$, with $\jsub{12}$.  Then this modified construction
will produce $256$ critical points below $\jsub{12}(\kappa_3) =
\gamma_{256}$, so these must be the first $256$ critical points.
But $\jsub{12} \Lequiv\jsub7(\kappa_2)/ j\jsub2$~\cite{\Dougherty,~\S5},
so the regular construction and the modified construction produce
exactly the same results below $\jsub7(\kappa_2) = \gamma_{256}$;
hence, the regular construction also produces the first
$256$ critical points. \QED

Since the irregular construction~\cite{\Dougherty,~\S4} actually has
$\jsub{12}$ as the first member of column~4, it produces the first
$256$ critical points as well as $\gamma_{256}$ itself.  It is possible
that the irregular construction produces a very large initial segment
of the list of critical points; if this is so, we would get equations
like $\jsub{11}(\gamma_3) = \gamma_{2^{2^{11}}}$ and
$\jsub{11}(\gamma_4) = \gamma_{2^{2^{2059}}}$.

\head 9. Conclusion \endhead
The results of this paper give some hope that an optimal recursive
construction producing all critical points is possible.  If
the patterns shown so far continue, then one can expect column
number~$N$ to contain embeddings which `look like' $\jsub{f(n)-N}$,
where $f(n)$ is a power of~2; this gives some guidance as to
the critical points patterns these columns should satisfy.

However, even if such a construction is completed, one will not
get a computation of the critical point counting function~$F$
unless one can tell where the critical points $\kappa_n$ occur
in the construction.  Since the function $F$ grows not only faster
than the Ackermann function, but faster than a slow iterate of
the Ackermann function~\cite{\Dougherty}, it must be the case
that the number of columns between the places where the
critical points $\kappa_n$ occur is also a fast-growing function.
This is a sign that there are still more complications to occur
later in the construction.